GeoClaw

\documentclass[pdflatex,sn-mathphys,iicol]{sn-jnl}

\usepackage[utf8]{inputenc}
\usepackage{graphicx} 

\usepackage{tikz,pgfplots}
\usepackage{program}
\usepackage{subcaption}
\usetikzlibrary{arrows, arrows.meta}
\usetikzlibrary{patterns}
\usepackage{geometry}
\usepackage{amssymb,amsmath}
\DeclareGraphicsExtensions{.pdf,.eps,.png,.jpg,.jpeg}
\usepackage[capitalize,nameinlink]{cleveref}

\jyear{2021}%

\theoremstyle{thmstyleone}%
\theoremstyle{thmstyletwo}%

\theoremstyle{thmstylethree}%

\begin{document}

\title[SRD Barrier Simulations]{Barrier Simulations with Shallow Water Equations using the State Redistribution Method}

\author*[1]{\fnm{Chanyang} \sur{Ryoo}}\email{cr2940@columbia.edu}

\author*[1]{\fnm{Kyle} \sur{Mandli}}\email{kyle.mandli@columbia.edu}


\affil*[1]{\orgdiv{APAM}, \orgname{Columbia University}, \orgaddress{\street{500 W 120th St}, \city{New York}, \postcode{10027}, \state{NY}, \country{USA}}}




\abstract{The representation of small scale barriers, such as sea-walls, in coastal flooding  simulations is a common computational constraint that can be difficult to overcome due to the combination of the need for resolution and the CFL constraining time stepping involved in solving the underlying PDE. This article proposes an approach that uses the state redistribution method (SRD) on shallow water equations to alleviate these problems while remaining conservative and geometrically easy to implement. We will demonstrate the method versus two other comparable approaches along with convergence results. In particular, we save computational time by reducing time steps fivefold when compared to simulation with same resolution using adaptive refinement at the barrier.}




\maketitle

\section{Introduction}\label{sec1}

The threat of coastal flooding to global coastlines is of growing concern.  Protecting those coastlines is increasingly important and simulating protective barriers at the coastlines at scale is becoming more critical.  The aim of this article is to provide improved computational methods for evaluating the effectiveness of these structures such as in \cite{usarmy,jmse8090725}. These structures tend to be narrow in one direction but long in another (e.g. the Maeslant barrier). The difficulty then comes down to how to deal with them computationally.  One approach is to use an unstructured grid with lot of refinement near the barrier. This poses two issues: first, the complexity of fitting a grid to various barrier shapes; second, the high computational cost spent resolving the thin barrier. Instead we propose a method that can accommodate a relatively unconstrained barrier across a quadrilateral grid while maintaining conservation properties and not being constrained by time stepping constraints that otherwise would be enforced by the CFL condition when resolving the barrier.

\begin{figure}[h!]
\centering
\begin{tikzpicture}[scale=0.7]
\draw[step=1cm,gray,very thin] (-1.5,-1.5) grid (3.5,3.5);
\draw[fill=blue!10] (0,0)--(0,0.48)--(-0.8,0.0)--(0,0);
\draw[fill=red!10] (0,0.48)--(-0.8,0.0)--(-1,0)--(-1,1)--(0,1)--(0,0.48);
\draw[fill=blue!10] (0,0)--(0,0.48)--(0.87,1.0)--(1,1)--(1,0)--(0,0);
\draw[fill=red!10] (0,0.48)--(0,1)--(0.87,1.0)--(0,0.48);
\draw[fill=blue!10] (0.87,1.0)--(1,1)--(1,1.08)--(0.87,1.0);
\draw[fill=red!10] (0,1)--(0,2)--(1,2)--(1,1.08)--(0.87,1.0)--(0,1);
\draw[fill=blue!10] (1,1)--(1,1.08)--(2,1.69)--(2,1)--(1,1);
\draw[fill=red!10] (1,1.08)-- (1,2)--(2,2)--(2,1.69)--(1,1.08);
\draw[fill=blue!10] (2,1)--(2,1.69)--(2.52,2.0)--(3,2)--(3,1)--(2,1);
\draw[fill=red!10] (2,2)--(2,1.69)--(2.52,2.0)--(2,2);
\draw[yellow,thick] (-1.3,-.3) -- (3.5,2.6);
\end{tikzpicture}
\caption{Some cut cells with two states (red and blue) generated by a linear barrier element.}
\label{fig:cut_cells}
\end{figure}
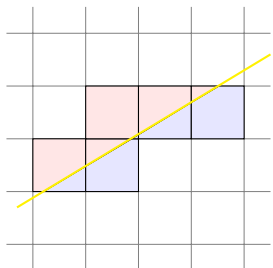

The shallow water equations in 2D are commonly used to model coastal flooding scenarios such as storm surge and tsunamis \cite{BERGER20111195,martin2010lake}. In this paper, instead of resolving a physical barrier, the surge barrier is modeled as a zero-width embedded boundary element that acts as a flux interface (line in \cref{fig:cut_cells}), similar to that done in \cite{zhang2013transition}. The interface mimics a solid barrier that either completely blocks off incoming waves or admits some waves over the barrier while reflecting others depending on its height. Our results show that this approximation are virtually the same as resolving the barrier, with added computational benefits.

Since we will be simulating our barriers within a Cartesian coordinate framework, the main problem that arises from this approximation is the \emph{small} or \emph{cut cell problem}, introduced by the barrier that cuts through arbitrary grid cells (\cref{fig:cut_cells}), which can get arbitraily small. This causes strict time step restriction via CFL condition. This stiffness problem can be solved by using a family of methods called \emph{cut cell methods} \cite{berger2012simplified,2006-Chung-p607,colella2006cartesian}. The idea is to apply special treatment to cut cells and their neighboring cells to lift the CFL restriction. With cut cell methods, there is no need to generate complex meshes that fit the barrier and allows for overall easier adaptation to different types of barriers \cite{BERGER20171,INGRAM2003561}. Cut cell methods therefore do not suffer from the CFL restriction that unstructured grids can have.

Here we present a simple finite volume method that uses the state redistribution cut cell method \cite{BERGER2020109820} and tackles piecewise linear barriers. Results show that even with considerably small cells, we can take full time steps. We treat two barrier problems: an angled barrier and a V-shaped barrier that has an obtuse angle at the corner. In \cref{sec1} we describe the model equation and problem. In \cref{sec2}, we explain the numerical method used to solve the problem in regular (uncut) cells and in \cref{sec3}, the state redistribution method (SRD) used in cut cells. In \cref{sec4} we present numerical examples and analysis of (1) a $20^\circ$ angled barrier and (2) a V-shaped barrier with an angle of $117^\circ$. Finally we present the computational efficiency that we achieve using our proposed method compared to the common method of refining the barrier in \cref{sec5}.

\section{Model equation and problem setup} \label{sec1}
\subsection{2D shallow water equations}
Shallow water equations (SWE) are a system of hyperbolic partial differential equations modeling horizontal, surface water motion. There are source terms to construct a geophysical version of SWE \cite{calhoun2008logically}, but for simplicity we only consider bathymetric source terms. SWE can be written as a set of conservation equations \cite{GDavid}:
\begin{align}
& h_t + (hu)_x + (hv)_y = 0 \\
& (hu)_t + \left(\frac{1}{2}gh^2 + hu^2 \right)_x + (huv)_y = -ghb_x\\
& (hv)_t + (huv)_x + \left(\frac{1}{2}gh^2+hv^2 \right)_y = -ghb_y,
\end{align}
for $(x,y) \in \Omega \subset \mathbb{R}^2$, where $h$ represents the height of water, $u$ and $v$ the $x$-velocity and $y$-velocity, $g$ the gravitational constant, and $b=b(x,y)$ the bathymetry. In the first equation, we have conservation of mass $h$ and in the second and third equations, if $\nabla b \equiv 0$, conservation of momentum. In the presence of bathymetric variation there is loss of conservation of momentum.

We can thus represent the SWE in the following form \cref{conseqn} as derived in \cite{leveque2002finite}. If $q=[h,hu,hv]$, $f(q)=[q_2, \frac{1}{2}gq_1^2 + \frac{q_2^2}{q_1},\frac{q_2q_3}{q_1}]$, $g(q)=[q_3,\frac{q_2q_3}{q_1},\frac{1}{2}gq_1^2 + \frac{q_3^2}{q_1}]$, and $\Psi(q,b)= [0,-gq_1b_x, -gq_1b_y]$, then we have for our SWE:
\begin{align}
    q_t + f(q)_x + g(q)_y = \Psi(q,b).
    \label{conseqn}
\end{align}

\subsection{Barrier representation on Cartesian grid}
\begin{figure}[h!]
\begin{subfigure}[b]{0.5\textwidth}
\centering
\begin{tikzpicture}
\draw[step=0.4cm] (0,0) grid (4,4);
\node at (-0.1, 1.2) {$1$};
\draw[red,thick] (0,1.2)--(4,2.5);
\node at (4.1, 2.5) {$2$};
\end{tikzpicture}
\caption{Model problem 1: slanted barrier.}
\label{fig:2d_setup1}
\end{subfigure}
 \hfill
\begin{subfigure}[b]{0.5\textwidth}
\centering
\begin{tikzpicture}
\draw[step=0.4cm] (0,0) grid (4,4);
\node at (-0.1, 2.7) {$1$};
\draw[red,thick] (0,2.7)--(2,1.27)--(4,2.7);
\node at (2, 1.2) {$2$};
\node at (4.1, 2.7) {$3$};

\end{tikzpicture}
\caption{Model problem 2: $V$-barrier}
\label{fig:2d_setup2}
\end{subfigure}
\caption{Grid setup of two model barrier problems. $1$, $2$, and $3$ denote the vertices of the barriers.}
\label{fig:2d_setup}
\end{figure}
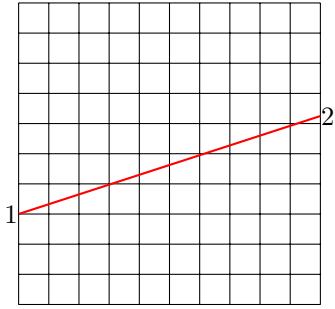
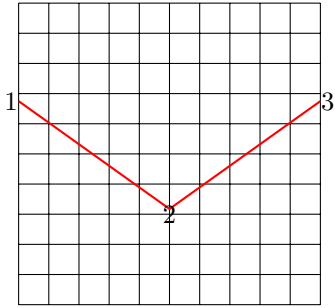

We now present the two model problems by showing their grid setup in \cref{fig:2d_setup}. Each grid cell has four volume averages: $H$, $HU$, $HV$ and $B$. We set a uniform bathymetry level $B=-2$ throughout the grid.

In red we have our zero width barriers. Note that the cut produces two states on either side of the line in the same cell position. This has both research and coding implications. It has research implications as this is the only instance known by the authors where there are state cells on both sides of a cut cell in 2D. In all cases that deal with cut cells \cite{CAUSON2000545,TUCKER2000591,may2017explicit}, only one side contains the state values. This also has coding implications as we must keep two arrays to track the upper lower state values in the same grid cell position.

\subsubsection{Computational aspects}
There are several input parameters for our model problem. First, we have the vertices of the barriers $1$, $2$, and $3$ (\cref{fig:2d_setup}) which determine the location and geometry of the linear and $V$-shaped barrier. Then we have the barrier height parameter $\beta$, which is the uniform height measured from the bathymetry level. For our model problems, we set vertex $1 = (0,0.3)$ and vertex $2=(0,0.653)$ for the linear barrier case ($20^\circ$) and vertex $1=(0,0.72)$, vertex $2=(0.5,0.412)$ and vertex $3=(1,0.72)$ for the V barrier case ($117^\circ$).

The barriers produce arbitrarily small cut cells on the order of 10$^{-5}\Delta x \Delta y$. We compute the locations of cut cells by using the vertices $A$, $B$, and $C$ of the barrier segments and the grid mesh size to identify intersections between the barrier and the grid. We compute the area of a cut cell $\mathcal{A}$ by then identifying its vertices $\{x_i,y_i\}_{i=1}^n$ and using the Shoelace formula. We pre-compute and save these for use in simulation runtime.
\section{Numerical method on uncut cells: wave propagation} \label{sec2}
We employ a finite volume method to solve SWE. The usual finite volume method is to divide the domain into grid cells, compute fluxes on each edge of the grid cell, and update the state via those fluxes.

Let $Q_{i,j}$ denote the state variable $[H_{i,j},(HU)_{i,j},(HV)_{i,j}]$, $F_{i\pm 1/2,j}$ denote the flux at the left/right edge of the cell, and $G_{i,j\pm 1/2}$ denote the flux at the bottom/upper edge (\cref{fig:FVM}). Then we have for a simple Godunov update:
\begin{align}
    \nonumber Q^{n+1}_{i,j} = Q^{n}_{i,j} - \frac{\Delta t}{\Delta x}(F^n_{i+1/2,j}-F^n_{i-1/2,j}) \\  - \frac{\Delta t}{\Delta y}(G^n_{i,j+1/2} - G^n_{i,j-1/2}).
\end{align}
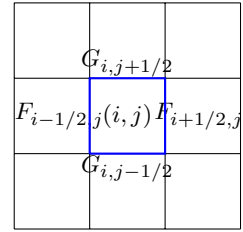
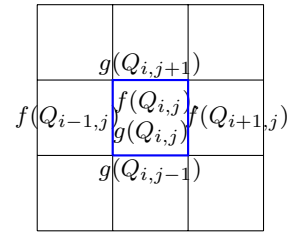
\begin{figure}[h!]
\begin{subfigure}[b]{0.5\textwidth}
\centering
\begin{tikzpicture}
\draw[step=1cm] (0,0) grid (6/2,6/2);
\node at (3/2, 3/2) {\scriptsize $(i,j)$};
\node at (1.4/2-0.1,3/2) {\scriptsize$F_{i-1/2,j}$};
\node at (4.7/2+0.1,3/2) {\scriptsize$F_{i+1/2,j}$};
\node at (3/2,1.6/2) {\scriptsize$G_{i,j-1/2}$};
\node at (3/2,4.4/2) {\scriptsize$G_{i,j+1/2}$};
\draw[blue,thick] (2/2,2/2)--(4/2,2/2)--(4/2,4/2)--(2/2,4/2)--(2/2,2/2);
\end{tikzpicture}
\caption{Grid edges of cell $(i,j)$ and its flux notation}
\label{fig:FVM}
\end{subfigure}
\begin{subfigure}[b]{0.5\textwidth}
\centering
\begin{tikzpicture}
\draw[step=1cm] (0,0) grid (6/2,6/2);
\node at (3/2+0.08, 3.2/2+0.1) {\scriptsize$f(Q_{i,j})$,};
\node at (3/2+0.01, 2.8/2-0.1) {\scriptsize$g(Q_{i,j})$};
\node at (1.2/2-0.2,3/2) {\scriptsize$f(Q_{i-1,j})$};
\node at (4.9/2+0.2,3/2) {\scriptsize$f(Q_{i+1,j})$};
\node at (3/2,1.6/2) {\scriptsize$g(Q_{i,j-1})$};
\node at (3/2,4.4/2) {\scriptsize$g(Q_{i,j+1})$};
\draw[blue,thick] (2/2,2/2)--(4/2,2/2)--(4/2,4/2)--(2/2,4/2)--(2/2,2/2);
\end{tikzpicture}
\caption{Flux vectors used in wave propagation for cell $(i,j)$}
\label{fig:FVM2}
\end{subfigure}
\caption{Flux scheme versus wave propagation}
\end{figure}%

However, there is another method called wave propagation, or flux-based wave decomposition \cite{bale2003wave}, which instead of using fluxes at the edges of a cell, linearly decomposes the difference of flux vectors from either side of each edge (\cref{fig:FVM2}). The matrices used for the linear decomposition are given by
\begin{align}
    A = \begin{bmatrix}
    1 & 0 & 1 \\
    U^*-\sqrt{gH^*} & 0 & U^* + \sqrt{gH^*} \\
    V^* & 1 & V^*
    \end{bmatrix}, \nonumber \\
    B = \begin{bmatrix}
    1 & 0 & 1 \\
    U^* & 1 & U^*  \\
    V^*-\sqrt{gH^*} & 0 & V^*+ \sqrt{gH^*}
    \end{bmatrix},
    \label{matrices}
\end{align}
where matrix $A$ is used to decompose difference of the flux vectors $f(q)$, and $B$ is used for $g(q)$, and $H^*,U^*,V^*$ represent special averages between $Q_{i,j}$ and neighbors, such as Roe or Einfeldt average \cite{einfeldt1988godunov,roe1981approximate}. Then the decomposition equations for the right and top edge become
\begin{align}
    f(Q_{i+1,j})-f(Q_{i,j})=A\gamma \label{diffcoeff}\\
    g(Q_{i,j+1})-g(Q_{i,j})=B\delta,\label{diffcoeff2}
\end{align}
where $\gamma, \, \delta \in \mathbb{R}^3.$ The reason for using the wave propagation method is that it is a versatile method that can also be applied to nonconservative hyperbolic equations \cite{leveque2002finite} and also an easily implemented method as long as the decomposing matrices are known and well-behaved.
The columns of $A$ and $B$ \cref{matrices} are the eigenvectors of the Jacobians $f'(q)$ and $g'(q)$, respectively. The eigenvalues corresponding to the $p^{\text{th}}$ column vector of $A$ are $\{(\sigma_A)_p\}_{p=1}^{3} =\{U^*-\sqrt{gH^*}, U^*, U^*+\sqrt{gH^*}\}$, and those corresponding to the $p^{\text{th}}$ column vector of $B$ are $\{(\sigma_B)_p\}_{p=1}^3 =\{V^*-\sqrt{gH^*},V^*, V^*+\sqrt{gH^*}\}$.

\subsection{Update formula using wave propagation}
Let $\{(\sigma_A)_p, A_{(\sigma_A)_p}\}$ and $\{(\sigma_B)_p,B_{(\sigma_B)_p}\}$ represent the $p^{\text{th}}$ eigenvalue-eigenvector pair of the Jacobians $f'(q)$ and $g'(q)$. Then we have the \emph{$A$-left} and \emph{$A$-right fluctuation waves} defined by
\begin{align}
    A^-\Delta Q_{i+1/2,j}=\sum_{p : (\sigma_A)_p <0} \gamma_p (\sigma_A)_p  A_{(\sigma_A)_p}\\
    A^+\Delta Q_{i-1/2,j} = \sum_{p : (\sigma_A)_p >0} \gamma_p (\sigma_A)_p A_{(\sigma_A)_p},
\end{align}
where $\gamma$ is the difference coefficient vector in \cref{diffcoeff}.

The \emph{$B$-up} and \emph{$B$-down fluctuation waves} are given by:
\begin{align}
B^-\Delta Q_{i,j-1/2} = \sum_{p : (\sigma_B)_p < 0} \delta_p (\sigma_B)_p B_{(\sigma_B)_p} \\
B^+\Delta Q_{i,j+1/2} = \sum_{p: (\sigma_B)_p >0} \delta_p (\sigma_B)_p B_{(\sigma_B)_p},
\end{align}
where $\delta$ is the difference coefficient vector in \cref{diffcoeff2}.
Then the state update formula via the wave propagation method becomes:
\begin{align}
\label{waveupd}
    Q^{n+1}_{i,j} = Q^n_{i,j} - \frac{\Delta t}{\Delta x}(A^- \Delta Q_{i+1/2,j} + A^+ \Delta Q_{i-1/2,j}) \nonumber \\ - \frac{\Delta t}{\Delta y}(B^- \Delta Q_{i,j+1/2} + B^+ \Delta Q_{i,j-1/2}).
\end{align}

\section{Numerical method on cut cells} \label{sec3}

\subsection{Attempt with $h$-boxes: Motivation for Using State Redistribution (SRD)}
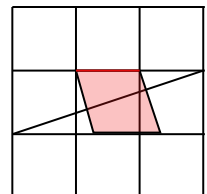
\begin{figure}[h!]
    \centering
    \tikzset{every picture/.style={line width=0.75pt}} 
\begin{tikzpicture}[scale=0.8,x=0.75pt,y=0.75pt,yscale=-1,xscale=1]

\draw  [draw opacity=0] (100,110) -- (221,110) -- (221,231) -- (100,231) -- cycle ; \draw   (100,110) -- (100,231)(140,110) -- (140,231)(180,110) -- (180,231)(220,110) -- (220,231) ; \draw   (100,110) -- (221,110)(100,150) -- (221,150)(100,190) -- (221,190)(100,230) -- (221,230) ; \draw    ;
\draw    (100,190) -- (220,150) ;
\draw  [fill={rgb, 255:red, 243; green, 27; blue, 27 }  ,fill opacity=0.27 ] (180,150) -- (193,189) -- (151,189) -- (140,150) -- cycle ;
\draw [color={rgb, 255:red, 224; green, 8; blue, 8 }  ,draw opacity=1 ]   (180,150) -- (140,150) ;

\end{tikzpicture}

    \caption{Ambiguous case of $h$-box averaging, required for edge in red.}
    \label{fig:2DV_hbox}
\end{figure}

The authors of this paper have attempted to solve this model problem using $h$-box type methods, with parallel, and diagonal barriers to deal with the small cell problem.
The difficulties of the $h$-box method applied to our model problem are (1) cumbersome geometrical calculations, (2) complicated update formulas, and (3) ambiguity in more complex cases. The biggest challenge of the $h$-box method is finding fluxes at the non-barrier edges of the cut cells (\cref{fig:2DV_hbox}). This is because $h$-box extensions off those edges will necessarily cross over the barrier, making it difficult to calculate their average. We will show how SRD overcomes all these challenges and can simulate arbitrary angled barrier and the $V$-shaped barrier.

\subsection{State Redistribution (SRD)}
\begin{figure}[h!]
 \tikzset{every picture/.style={line width=0.75pt}} 
\centering
\begin{tikzpicture}[scale=0.8,x=0.75pt,y=0.75pt,yscale=-1,xscale=1]

\draw   (338.2,111.1) -- (388.2,111.1) -- (388.2,161.1) -- (338.2,161.1) -- cycle ;
\draw    (171,129) -- (193.2,112.6) ;
\draw   (256,112) -- (306,112) -- (306,162) -- (256,162) -- cycle ;
\draw    (269.2,161.6) -- (288.2,111.6) ;
\draw   (170,112) -- (220,112) -- (220,162) -- (170,162) -- cycle ;
\draw    (338.2,143.6) -- (388.2,126.6) ;
\draw   (427,110) -- (477,110) -- (477,160) -- (427,160) -- cycle ;
\draw    (456.2,160.6) -- (476.2,135.6) ;
\draw   (170,181) -- (220,181) -- (220,231) -- (170,231) -- cycle ;
\draw    (198.6,181.6) -- (221,195) ;
\draw   (256,181) -- (306,181) -- (306,231) -- (256,231) -- cycle ;
\draw    (271.6,181.6) -- (291.6,231.6) ;
\draw   (338,180) -- (388,180) -- (388,230) -- (338,230) -- cycle ;
\draw    (338.6,196.6) -- (387.6,208.6) ;
\draw   (427,180) -- (477,180) -- (477,230) -- (427,230) -- cycle ;
\draw    (427.6,216.6) -- (446,230) ;
\end{tikzpicture}
\caption{All possible types of cut cells for both upper and lower cut cells for a barrier segment \cite{causon2000calculation}.}
\label{fig:cutcells_ex}
\end{figure}
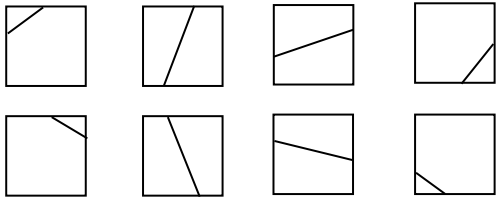

The equation \cref{waveupd} shown in the previous section does not work on cut cells because of the geometry of the cut cells. We have different kinds of cut cells as shown in \cref{fig:cutcells_ex}.

Instead of \cref{waveupd}, we apply the state redistribution method. The main idea of the method is to first do a conservative but unstable update on all cut cells and correct for the instability by doing a stable redistribution of those updates.
\subsection{1D SRD}
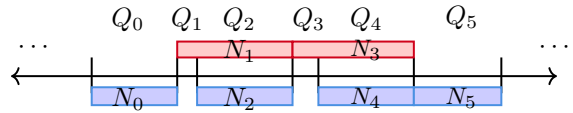
\begin{figure}[h!]
    \centering
\tikzset{every picture/.style={line width=0.75pt}} 
\begin{tikzpicture}[x=0.75pt,y=0.75pt,yscale=-1,xscale=1]

\draw  [<->]  (190,104) -- (464.6,104) ;
\draw    (230.6,112.6) -- (230.6,94.6) ;
\draw    (283.6,112.6) -- (283.6,94.6) ;
\draw    (331.6,112.6) -- (331.6,94.6) ;
\draw    (344.6,112.6) -- (344.6,94.6) ;
\draw    (392.6,112.6) -- (392.6,94.6) ;
\draw    (273.6,112.6) -- (273.6,94.6) ;
\draw    (436.6,112.6) -- (436.6,94.6) ;
\draw  [color={rgb, 255:red, 74; green, 144; blue, 226 }  ,draw opacity=1 ,fill=blue!20] (230.6,109.6) -- (273.6,109.6) -- (273.6,118.6) -- (230.6,118.6) -- cycle ;
\draw  [color={rgb, 255:red, 74; green, 144; blue, 226 }  ,draw opacity=1 ,fill=blue!20] (283.6,109.6) -- (331.6,109.6) -- (331.6,118.6) -- (283.6,118.6) -- cycle ;
\draw  [color={rgb, 255:red, 74; green, 144; blue, 226 }  ,draw opacity=1 ,fill=blue!20] (344.6,109.6) -- (392.6,109.6) -- (392.6,118.6) -- (344.6,118.6) -- cycle ;
\draw  [color={rgb, 255:red, 74; green, 144; blue, 226 }  ,draw opacity=1 ,fill=blue!20] (392.6,109.6) -- (436.6,109.6) -- (436.6,118.6) -- (392.6,118.6) -- cycle ;
\draw  [color={rgb, 255:red, 208; green, 2; blue, 27 }  ,draw opacity=1 ,fill=red!20] (273.6,86.6) -- (331.6,86.6) -- (331.6,94.6) -- (273.6,94.6) -- cycle ;
\draw  [color={rgb, 255:red, 208; green, 2; blue, 27 }  ,draw opacity=1,fill=red!20 ] (331.6,86.6) -- (392.6,86.6) -- (392.6,94.6) -- (331.6,94.6) -- cycle ;


\draw (239,68.4) node [anchor=north west][inner sep=0.75pt]    {$Q_{0}$};
\draw (239,108.4) node [anchor=north west][inner sep=0.75pt]    {$N_{0}$};

\draw (269,68.4) node [anchor=north west][inner sep=0.75pt]    {$Q_{1}$};
\draw (330,68.4) node [anchor=north west][inner sep=0.75pt]    {$Q_{3}$};
\draw (295,68.4) node [anchor=north west][inner sep=0.75pt]    {$Q_{2}$};
\draw (295,84.4) node [anchor=north west][inner sep=0.75pt]    {$N_{1}$};
\draw (295,108.4) node [anchor=north west][inner sep=0.75pt]    {$N_{2}$};
\draw (359,68.4) node [anchor=north west][inner sep=0.75pt]    {$Q_{4}$};
\draw (359,84.4) node [anchor=north west][inner sep=0.75pt]   {\small$N_{3}$};
\draw (407,67.4) node [anchor=north west][inner sep=0.75pt]    {$Q_{5}$};
\draw (407,108.4) node [anchor=north west][inner sep=0.75pt]    {\small $N_{5}$};
\draw (359,108.4) node [anchor=north west][inner sep=0.75pt]    {\small $N_{4}$};
\draw (192,85) node [anchor=north west][inner sep=0.75pt]   [align=left] {$\cdots$};
\draw (454,85) node [anchor=north west][inner sep=0.75pt]   [align=left] {$\cdots$};
\end{tikzpicture}
    \caption{An irregular grid with small cells and the neighborhoods used in SRD in blue and red.}
    \label{fig:1DSRD}
\end{figure}

The main components of the method are best explained in 1D. In \cref{fig:1DSRD} we have an irregular grid with two small cells, $Q_1$ and $Q_3$, with areas $\alpha \Delta x$ ($\alpha<0.5$) \cite{BERGER2015180}, where a normal update cannot be taken due to the CFL condition. Every other grid cell has length $\Delta x$. This is the same grid setup as in the example in \cite{BERGER2020109820}, but we provide a slightly different method which will also apply to our 2D model problems.

Before applying the state redistribution, each cell, including the small cells, is given a regular update with the $\Delta t$ prescribed by CFL condition with a regular $\Delta x$, \emph{regardless} of cell size:
$$
Q^{n+1}_i = Q^{n}_i - \frac{\Delta t}{\alpha_i \Delta x}(A^+\Delta Q_{i-1/2} + A^-\Delta Q_{i+1/2}),
$$
where $\alpha_i=\alpha$ when $i=1,3$ and $\alpha_i=1$ when $i\ne 1,3.$
This will be an unstable yet conservative update for the small cells.

Then we define \emph{neighborhoods} for each cell, $\{N_i\}$, which are collection of cells to achieve a collective area greater than $0.5 \Delta x$. This means that for cells with $\alpha \ge 0.5$, they are their own neighborhood, e.g. $N_{-1} = Q_{-1}$ (blue in \cref{fig:1DSRD}). For the small cells, we can select the regular sized cell immediately to their right and this achieves a collective area larger than $0.5 \Delta x$. The two red boxes in \cref{fig:1DSRD} are such neighborhoods for $Q_1$ and $Q_3$. This is different from the example shown in \cite{BERGER2020109820} where neighboring cells are taken from both sides. However, as we shall see, forming this one sided neighborhood is conservative and will be needed in our model problem due to the presence of the zero width barrier.

The neighborhoods also have a special average. The average of a neighborhood is calculated by using its comprising cells' averages, areas, and their \emph{overlap counts}, which is the total number of neighborhoods lying over them. For instance, cells $Q_2$ and $Q_4$ both have overlap count of $2$ ($N_1,N_2$ and $N_3,
N_4$ in \cref{fig:1DSRD}). This gives the following neighborhood averages:
\begin{align}
    &N_1 = \frac{1}{\alpha + \frac{1}{2}}(\alpha Q^{n+1}_1 + \frac{1}{2} Q^{n+1}_2), \\
    &N_3 = \frac{1}{\alpha + \frac{1}{2}}(\alpha Q^{n+1}_3 + \frac{1}{2} Q^{n+1}_4),\\
    &N_i = Q^{n+1}_i \, \text{for $i\ne 1,3$}.
\end{align}
Note that $\Delta x$ factors out and weights for $Q_2,Q_4$ have overlap count 2 in the denominator. Discounting by overlap count gives the effective volume and is used for conservation purposes as can be seen by checking the mass from before to after SRD stabilization:
$$
\hspace{-2cm}     \overbrace{Q_0 + \alpha Q_1 + Q_2 + \alpha Q_3 + Q_4 + Q_5}^{\text{pre-SRD}}
$$
$$
    =Q_{0} +\underbrace{ \alpha N_1 + \frac{1}{2}(N_1+N_2) + \alpha N_3 + \frac{1}{2}(N_3+N_4) }_{SRD}+Q_5
$$
$$
    =Q_0 + (\alpha Q_1 + \frac{1}{2}Q_2) + \frac{1}{2} Q_2+ (\alpha Q_3 + \frac{1}{2}Q_4) + \frac{1}{2} Q_4 + Q_5
$$\vspace{0.1cm}
$$
    = \underbrace{Q_0 + \alpha Q_1 + Q_2 + \alpha Q_3 + Q_4 + Q_5.}_{\text{post-SRD}}
$$

Finally, stabilized updates of the cells are found by using the neighborhood average values lying over the cells and averaging them if there are multiple neighborhoods:
\begin{align}
    &Q_2^{n+1} = \frac{1}{2}(N_1+N_2) \\
    &Q_4^{n+1} = \frac{1}{2}(N_3+N_4) \\
    &Q_i^{n+1} = N_i \, \text{for $i \ne 2,4$}.
\end{align}

\subsection{The original 2D SRD method}
In this section we describe an adapted version of original 2D SRD method applicable for impermeable solid boundaries. The method here will apply to our case with the permeable barrier interface, with the exception of the barrier edge handling.

The added complexity is now finding the conservative but unstable updates of the geometrically specific cut cells during the first step of SRD method. There are 8 total types of cut cells (4 if we disregard the sign of the cut slope) as seen in \cref{fig:cutcells_ex}. This means that we must update the cut cells in at most 4 different ways depending on their shapes.

\subsubsection{Conservative but unstable update}
\begin{figure}[h!]
    \tikzset{every picture/.style={line width=0.75pt}} 
\centering
\begin{subfigure}[b]{0.4\textwidth}

\begin{tikzpicture}[scale=0.8,x=0.75pt,y=0.75pt,yscale=-1,xscale=1]
\hspace{-1cm}
\draw  [draw opacity=0] (231.6,81.6) -- (382.6,81.6) -- (382.6,232.6) -- (231.6,232.6) -- cycle ; \draw   (231.6,81.6) -- (231.6,232.6)(281.6,81.6) -- (281.6,232.6)(331.6,81.6) -- (331.6,232.6)(381.6,81.6) -- (381.6,232.6) ; \draw   (231.6,81.6) -- (382.6,81.6)(231.6,131.6) -- (382.6,131.6)(231.6,181.6) -- (382.6,181.6)(231.6,231.6) -- (382.6,231.6) ; \draw    ;
\draw[fill=gray!20]   (365.6,81.9) -- (231.6,199.6)--(231.6,81.9) ;
\draw    (189.6,168.6) .. controls (215.93,191.96) and (237.17,203.37) .. (276.41,174.49) ;
\draw [shift={(277.6,173.6)}, rotate = 503.13] [color={rgb, 255:red, 0; green, 0; blue, 0 }  ][line width=0.75]    (10.93,-3.29) .. controls (6.95,-1.4) and (3.31,-0.3) .. (0,0) .. controls (3.31,0.3) and (6.95,1.4) .. (10.93,3.29)   ;
\draw    (415.6,102.6) .. controls (372.04,50.13) and (313.78,80.97) .. (324.26,129.14) ;
\draw [shift={(324.6,130.6)}, rotate = 256.24] [color={rgb, 255:red, 0; green, 0; blue, 0 }  ][line width=0.75]    (10.93,-3.29) .. controls (6.95,-1.4) and (3.31,-0.3) .. (0,0) .. controls (3.31,0.3) and (6.95,1.4) .. (10.93,3.29)   ;

\draw (122,130) node [anchor=north west][inner sep=0.75pt]   [align=left] {};
\draw (283,160.4) node [anchor=north west][inner sep=0.75pt]    {$l_{1}$};
\draw (316,133.4) node [anchor=north west][inner sep=0.75pt]    {$l_{2}$};
\draw (151,140.4) node [anchor=north west][inner sep=0.75pt]    {$l_{1} A^{+} \Delta Q_{i}{}_{-1}{}_{/}{}_{2}{}_{,}{}_{j}$};
\draw (373,106.4) node [anchor=north west][inner sep=0.75pt]    {$l_{2} B^{-} \Delta Q_{i}{}_{,}{}_{j+1/2}$};
\draw (393,146.4) node [anchor=north west][inner sep=0.75pt]    {$j$};
\draw (304,237.4) node [anchor=north west][inner sep=0.75pt]    {$i$};

\end{tikzpicture}

    \caption{The weighted waves used to update the cut cell $(i,j)$. The wave from the left edge is weighted by $l_1$ and the wave from the top edge is weighted by $l_2$.}
    \label{fig:sidelen}
    \end{subfigure}
    \centering
    \hfill
\hspace{-1cm}
\begin{subfigure}[b]{0.4\textwidth}
\centering

\begin{tikzpicture}[scale=0.8,x=0.75pt,y=0.75pt,yscale=-1,xscale=1]

\draw  [draw opacity=0] (231.6,81.6) -- (382.6,81.6) -- (382.6,232.6) -- (231.6,232.6) -- cycle ; \draw   (231.6,81.6) -- (231.6,232.6)(281.6,81.6) -- (281.6,232.6)(331.6,81.6) -- (331.6,232.6)(381.6,81.6) -- (381.6,232.6) ; \draw   (231.6,81.6) -- (382.6,81.6)(231.6,131.6) -- (382.6,131.6)(231.6,181.6) -- (382.6,181.6)(231.6,231.6) -- (382.6,231.6) ; \draw    ;
\draw[fill=gray!20]   (365.6,81.6) -- (231.6,199.6)--(231.6,81.6)  ;
\draw[fill=white]   (277.6,97.6) -- (308.6,131.6) -- (281.6,155.6) -- (250.6,121.6) -- cycle ;
\draw   (308.6,131.6) -- (339.6,165.6) -- (312.6,189.6) -- (281.6,155.6) -- cycle ;
\draw    (401.6,109.6) .. controls (378.83,64.06) and (292.35,96.93) .. (294.99,142.22) ;
\draw [shift={(295.1,143.6)}, rotate = 264.40999999999997] [color={rgb, 255:red, 0; green, 0; blue, 0 }  ][line width=0.75]    (10.93,-3.29) .. controls (6.95,-1.4) and (3.31,-0.3) .. (0,0) .. controls (3.31,0.3) and (6.95,1.4) .. (10.93,3.29)   ;

\draw (300.6,144) node [anchor=north west][inner sep=0.75pt]    {$l_{3}$};
\draw (393,146.4) node [anchor=north west][inner sep=0.75pt]    {$j$};
\draw (304,237.4) node [anchor=north west][inner sep=0.75pt]    {$i$};
\draw (385,110.4) node [anchor=north west][inner sep=0.75pt]    {$l_{3} A^{+} \Delta \tilde{Q}_{i,j}$};

\end{tikzpicture}
\caption{The normal averages used to produce waves at the barrier edge to update the cut cell $(i,j)$. The wave from the barrier edge is weighted by $l_3.$}
    \label{fig:barlen}

\end{subfigure}
\caption{Ratio of edge lengths to mesh size as weights on waves from edges of a cut cell.}
    \label{fig:cutcellupd}
\end{figure}
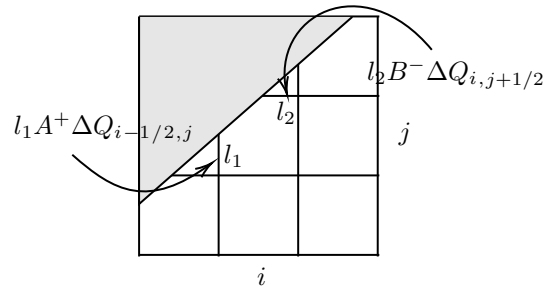
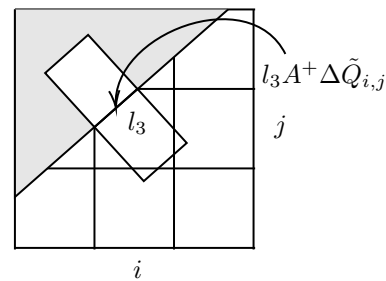

The way we compute the conservative but unstable updates is by first identifying the cut cell edge lengths. Using the side lengths, we weight the waves arising from those edges as seen in \cref{fig:cutcellupd}. The waves from the vertical/horizontal edges are computed by using wave propagation with the two adjacent cell values.

At the barrier edge, the original SRD method computes fluctuation waves by first creating normal ``$h$-boxes'' \cite{berger2012simplified} extending towards both sides of the barrier edge as in \cref{fig:barlen}. The states of these $h$-boxes are computed by volume weighting the averages of the covered cells and by rotating their momentum with respect to the barrier direction, denoted by $\tilde{Q}$ in \cref{fig:barlen}. The average of the $h$-box intruding into the solid region is computed by simply negating the velocity in the normal direction to the barrier (\cref{fig:barlen}). Since SWE is rotationally invariant, the waves at the barrier edge are calculated using method described in \cref{sec2} with these two $h$-box averages.

Once the cut cells are updated conservatively using the weighted waves at each edge, neighborhoods are found for each cell and compute the neighborhood averages for the SRD updates. The use of overlap counts and small cell criteria ($\mathcal{A} < 0.5 \Delta x \Delta y$) are all the same, so we only show how neighborhoods are formed in the 2D case in \cref{nhoodsec}.

\subsubsection{Finding neighborhoods}
\label{nhoodsec}
The two possible types of neighborhood are the \emph{normal} and the \emph{grid neighborhood}. As can be seen in left subfigure of \cref{fig:2Dnhood}, a normal neighborhood of cell $(i,j)$ only uses the cell directly below the cut cell. Note that the neighboring cell could also be a cut cell. As long as the total area of the neighborhood exceeds $0.5\Delta x \Delta y$, it is a valid neighborhood. There are cases when just including the direct neighboring cell does not result in a neighborhood of size bigger than $0.5\Delta x \Delta y,$ as in the right subfigure in \cref{fig:2Dnhood}. In this case, the neighborhood of cell $(i,j)$ needs not only cell $(i,j-1)$ but also the cells $(i-1,j-1)$ and $(i-1,j)$ to be included in their neighborhood. This type of neighborhood is used for cut cells that are especially constricted as in the figure and we will not use this type as our V barrier is obtuse, which makes physical sense.

\begin{figure}[h!]
    \centering
    \tikzset{every picture/.style={line width=0.75pt}} 
\begin{subfigure}[b]{0.4\textwidth}
\centering
\begin{tikzpicture}[scale=0.7,x=0.75pt,y=0.75pt,yscale=-1,xscale=1]

\draw (61,147.4) node [anchor=north west][inner sep=0.75pt]    {$j$};
\draw (151,237.4) node [anchor=north west][inner sep=0.75pt]    {$i$};
\draw  [draw opacity=0] (81.6,81.6) -- (232.6,81.6) -- (232.6,232.6) -- (81.6,232.6) -- cycle ; \draw   (81.6,81.6) -- (81.6,232.6)(131.6,81.6) -- (131.6,232.6)(181.6,81.6) -- (181.6,232.6)(231.6,81.6) -- (231.6,232.6) ; \draw   (81.6,81.6) -- (232.6,81.6)(81.6,131.6) -- (232.6,131.6)(81.6,181.6) -- (232.6,181.6)(81.6,231.6) -- (232.6,231.6) ; \draw    ;
\draw [fill=gray!20,opacity=0.8]   (231.6,109.6) -- (93.6,231.6)--(81.6,232.6)--(81.6,81.6)--(231.6,81.6) ;
\draw [green,fill=green!10,opacity=0.8]  (181.6,153.6) -- (181.6,231.6) -- (131.6,231.6) -- (131.6,197.6) -- cycle ;
\end{tikzpicture}
\end{subfigure}
\begin{subfigure}[b]{0.4\textwidth}
\begin{tikzpicture}[scale=0.7,x=0.75pt,y=0.75pt,yscale=-1,xscale=1]

  \draw  [draw opacity=0] (326.6,80.6) -- (627.6,80.6) -- (627.6,231.6) -- (326.6,231.6) -- cycle ; \draw   (326.6,80.6) -- (326.6,231.6)(376.6,80.6) -- (376.6,231.6)(426.6,80.6) -- (426.6,231.6)(476.6,80.6) -- (476.6,231.6)(526.6,80.6) -- (526.6,231.6)(576.6,80.6) -- (576.6,231.6)(626.6,80.6) -- (626.6,231.6) ; \draw   (326.6,80.6) -- (627.6,80.6)(326.6,130.6) -- (627.6,130.6)(326.6,180.6) -- (627.6,180.6)(326.6,230.6) -- (627.6,230.6) ; \draw    ;
\draw    (439.6,80.6) -- (476.6,230.6) ;
\draw [green, fill=green!10,opacity=0.8, thick](449.6,130.6) -- (498.6,130.6)-- (476.6,230.6);
\draw    (509.6,80.6) -- (476.6,230.6) ;
\draw[fill=gray!20,opacity=0.8] (439.6,80.6) -- (476.6,230.6)--(326.6,231.6)  -- (326.6,80.6)  ;
\draw[fill=gray!20,opacity=0.8] (509.6,80.6) -- (476.6,230.6) --(627.6,230.6) --(627.6,80.6)   -- (376.6,80.6)  ;
\draw (326.6,80.6)--(409.6,80.6);

\draw (303,197.4) node [anchor=north west][inner sep=0.75pt]    {$j$};
\draw (497,236.4) node [anchor=north west][inner sep=0.75pt]    {$i$};
\draw (437,237.4) node [anchor=north west][inner sep=0.75pt]    {$i-1$};
\end{tikzpicture}

\end{subfigure}
    \caption{Two types of neighborhoods in 2D SRD. The solid region is made opaque for clarity in identifying the indices.}
    \label{fig:2Dnhood}
\end{figure}

\subsection{SRD applied to model problem}
In our model problem the main adjustment we need to make to the original SRD method is handling the presence of the additional state on the other side of the barrier. We denote the upper cut cell by the superscript $Q^U_{i,j}$ and the lower cut cell by the superscript $Q^L_{i,j}.$ As shown in \cref{fig:dblstate}, we remark that aside from the barrier edge, the weighted wave propagation method as illustrated in \cref{fig:sidelen} applies to both the upper cut cell and the lower cut cell at each edge.

We now need to allow communication between the two sides of the cut through the barrier edge. Consequently, the key differences from the original SRD method are (1) the states used for the Riemann problem at the barrier edge and (2) calculation of fluctuation at the barrier edge, called \emph{wave redistribution}. Wave redistribution was first developed in \cite{li2021h}, but we develop it for the first time in 2D.
\begin{figure}[h!]
    \tikzset{every picture/.style={line width=0.75pt}} 
\centering
\begin{subfigure}[b]{0.4\textwidth}

\begin{tikzpicture}[scale=0.8,x=0.75pt,y=0.75pt,yscale=-1,xscale=1]

\draw  [draw opacity=0] (241.6,74.6) -- (392.6,74.6) -- (392.6,225.6) -- (241.6,225.6) -- cycle ; \draw   (241.6,74.6) -- (241.6,225.6)(291.6,74.6) -- (291.6,225.6)(341.6,74.6) -- (341.6,225.6)(391.6,74.6) -- (391.6,225.6) ; \draw   (241.6,74.6) -- (392.6,74.6)(241.6,124.6) -- (392.6,124.6)(241.6,174.6) -- (392.6,174.6)(241.6,224.6) -- (392.6,224.6) ; \draw    ;
\draw    (379.6,74.6) -- (241.6,196.6) ;
\draw    (264.6,244.4) .. controls (246.78,202.82) and (247.58,193.58) .. (290.68,162.94) ;
\draw [shift={(292,162)}, rotate = 504.73] [color={rgb, 255:red, 0; green, 0; blue, 0 }  ][line width=0.75]    (10.93,-3.29) .. controls (6.95,-1.4) and (3.31,-0.3) .. (0,0) .. controls (3.31,0.3) and (6.95,1.4) .. (10.93,3.29)   ;
\draw    (407.6,170.4) .. controls (426.5,102.74) and (391.95,81.61) .. (333.88,123.37) ;
\draw [shift={(333,124)}, rotate = 323.98] [color={rgb, 255:red, 0; green, 0; blue, 0 }  ][line width=0.75]    (10.93,-3.29) .. controls (6.95,-1.4) and (3.31,-0.3) .. (0,0) .. controls (3.31,0.3) and (6.95,1.4) .. (10.93,3.29)   ;
\draw    (315.6,46.4) .. controls (255.82,-2.6) and (281.52,79.97) .. (298.57,104) ;
\draw [shift={(299.6,105.4)}, rotate = 232.31] [color={rgb, 255:red, 0; green, 0; blue, 0 }  ][line width=0.75]    (10.93,-3.29) .. controls (6.95,-1.4) and (3.31,-0.3) .. (0,0) .. controls (3.31,0.3) and (6.95,1.4) .. (10.93,3.29)   ;
\draw    (218.6,112.4) .. controls (219.59,122.3) and (222.15,169.84) .. (272.07,143.23) ;
\draw [shift={(273.6,142.4)}, rotate = 511.08] [color={rgb, 255:red, 0; green, 0; blue, 0 }  ][line width=0.75]    (10.93,-3.29) .. controls (6.95,-1.4) and (3.31,-0.3) .. (0,0) .. controls (3.31,0.3) and (6.95,1.4) .. (10.93,3.29)   ;

\draw (396,140.4) node [anchor=north west][inner sep=0.75pt]    {$j$};
\draw (310,226.4) node [anchor=north west][inner sep=0.75pt]    {$i$};
\draw (295,152.4) node [anchor=north west][inner sep=0.75pt]    {$l_{1}$};
\draw (323,124.4) node [anchor=north west][inner sep=0.75pt]    {$l_{2}$};
\draw (277,128.4) node [anchor=north west][inner sep=0.75pt]    {$l_{3}$};
\draw (301.6,107.8) node [anchor=north west][inner sep=0.75pt]    {$l_{4}$};
\draw (230,254.4) node [anchor=north west][inner sep=0.75pt]    {$l_{1} A^{+} \Delta Q_{i-1/2,j}^{L}$};
\draw (375,173.4) node [anchor=north west][inner sep=0.75pt]    {$l_{2} B^{-} \Delta Q_{i,j+1/2}^{L}$};
\draw (296,50.4) node [anchor=north west][inner sep=0.75pt]    {$l_{4} B^{-} \Delta Q_{i,j+1/2}^{U}$};
\draw (183,92.4) node [anchor=north west][inner sep=0.75pt]    {$l_{3} A^{+} \Delta Q_{i-1/2,j}^{U}$};

\end{tikzpicture}

    \caption{The weighted waves used to update both cut cells $Q_{i,j}^U$ and $Q_{i,j}^L$ at the non-barrier edges. Weights $l_i$ are lengths of the cut cell edges.}
    \label{fig:dblstate}
    \end{subfigure}
    \centering
    \hfill
\begin{subfigure}[b]{0.4\textwidth}
\centering

\begin{tikzpicture}[scale=0.8,x=0.75pt,y=0.75pt,yscale=-1,xscale=1]

\draw  [draw opacity=0] (231.6,81.6) -- (382.6,81.6) -- (382.6,232.6) -- (231.6,232.6) -- cycle ; \draw   (231.6,81.6) -- (231.6,232.6)(281.6,81.6) -- (281.6,232.6)(331.6,81.6) -- (331.6,232.6)(381.6,81.6) -- (381.6,232.6) ; \draw   (231.6,81.6) -- (382.6,81.6)(231.6,131.6) -- (382.6,131.6)(231.6,181.6) -- (382.6,181.6)(231.6,231.6) -- (382.6,231.6) ; \draw    ;
\draw   (365.6,81.6) -- (231.6,199.6)  ;
\draw    (401.6,109.6) .. controls (378.83,64.06) and (292.35,96.93) .. (294.99,142.22) ;
\draw [shift={(295.1,143.6)}, rotate = 264.40999999999997] [color={rgb, 255:red, 0; green, 0; blue, 0 }  ][line width=0.75]    (10.93,-3.29) .. controls (6.95,-1.4) and (3.31,-0.3) .. (0,0) .. controls (3.31,0.3) and (6.95,1.4) .. (10.93,3.29)   ;

\draw (393,146.4) node [anchor=north west][inner sep=0.75pt]    {$j$};
\draw (304,237.4) node [anchor=north west][inner sep=0.75pt]    {$i$};
\draw (385,110.4) node [anchor=north west][inner sep=0.75pt]    {$l_{5} A^{\pm} \Delta \tilde{Q}_{i,j}$};
\draw (300.6,154) node [anchor=north west][inner sep=0.75pt]    {$\tilde{Q}^{L}_{i,j}$};
\draw (260.6,132) node [anchor=north west][inner sep=0.75pt]    {$\tilde{Q}^{U}_{i,j}$};

\end{tikzpicture}
\caption{The rotated averages $\tilde{Q}_{i,j}^{U/L}$ of both cut cells are used to produce waves at the barrier edge to update the cut cells $Q_{i,j}^{U/L}$. The wave is computed using wave redistribution and weighted by barrier edge length $l_5.$}
    \label{fig:barst}

\end{subfigure}
\caption{Weighted waves with edge lengths of a cut cell for both upper and lower cut cells.}
    \label{fig:dblstates}
\end{figure}

\subsubsection{Wave redistribution at barrier edge}
\label{WR}
For the Riemann problem at the barrier edge, we cannot simply negate the normal momentum. We must define a new Riemann problem to enable overtopping. First, we do not need to use the normal $h$-boxes to compute the fluctuation. Instead we can simply use the rotated state variables $\tilde{Q}^U_{i,j}, \tilde{Q}^L_{i,j}$ on either side of the barrier cut:
\begin{align}
    \tilde{Q}^h_{i,j} = R_{i,j} Q^h_{i,j},
\end{align}
where $h=L$ or $U,$ and
\begin{align*}
    R_{i,j} = \begin{bmatrix}
    1 & 0 & 0 \\
    0 & \hat{n}_1 & \hat{n}_2 \\
    0 & \hat{t}_1 & \hat{t}_2
    \end{bmatrix},
\end{align*}
where rotation vectors $\hat{n}, \hat{t}$ are simply orthonormal with respect to the barrier. This simple usage of state variables avoids any complicated geometrical computation at the V-tip.

Once we have rotated the states, these two states become the left and right states of two ghost problems, which are needed here just as in any boundary condition. The wave redistribution then redistributes waves from the two ghost Riemann problems. The ghost cell $Q^*$ introduced in between the two rotated states is a cell with $B^*=\beta$, the barrier height, and $Q^*=0$ if $\beta > \min (H_{i,j}^U,H_{i,j}^L)$, to mimic a dry state on top of the barrier. Otherwise, we let the ghost height $H^*=\min(H_{i,j}^U-\beta,H_{i,j}^L-\beta)$ (minimum overtopping amount) and $HU^*, \, HV^* = \min (HU^L_{i,j},HU^U_{i,j}), \, \min (HV^L_{i,j},HV^U_{i,j})$.

With the ghost cell, we have the following decompositions to perform:
\begin{align}
    f(\tilde{Q}^{L}_{i,j})-f(Q^*) - \Psi_L = A_L\gamma \\
    f(Q^*)-f(\tilde{Q}^{U}_{i,j}) - \Psi_U = A_U\delta,
\end{align}
where $A_L$ is the matrix similar to \cref{matrices} for the Riemann problem between $\tilde{Q}^{L}_{i,j}$ and $Q^*$ with associated eigenvalues and eigenvectors $\{\sigma_L^i, \mathbf{r}_L^i\}$ and $A_U$ is that for the Riemann problem between $\tilde{Q}^{U}_{i,j}$ and $Q^*$ with eigenvalues and eigenvectors $\{\sigma_U^i, \mathbf{r}_U^i\}$. The terms $\Psi_L$ and $\Psi_U$ are given by $\Psi_L = [0,g\overline{H}_L(B_L-B^*),0]$ and $\Psi_U=[0,g\overline{H}_U(B^*-B_U),0]$, where $\overline{H}_{L/U}$ is $0.5(H_{L/U}+H^*)$. This simple subtraction of source term from the flux differences is one of the key merits of the wave decomposition method \cite{bale2003wave}.

Then wave redistribution sets a new set of eigenvalues and eigenvectors $\{\omega_i, \boldsymbol{\rho}_i\}$ as follows:
\begin{align}
    &\omega_i = \frac{1}{2} ( \sigma_U^i + \sigma_L^i) \\
    &\boldsymbol{\rho}_i = [1,\omega_i,\frac{1}{2}(\tilde{V}^L + \tilde{V}^U)] \text{ for $i=1,3$} \\
    &\boldsymbol{\rho}_2 = [0,0,1],
\end{align}
where $\tilde{V}^{L/U}$ is the rotated transverse velocity of $\tilde{Q}^{L/U}_{i,j}.$

Finally, wave redistribution solves
\begin{align}
    [A_L \| A_U] (\gamma : \delta) = \boldsymbol{\mathrm{P}} \epsilon,
\end{align}
where $[A_L \| A_U]$ is the augmented matrix of two matrices $A_L$ and $A_U$, $\gamma : \delta$ is the augmented vector of coefficient vectors $\gamma$ and $\delta$, along the same column axis, and $\boldsymbol{\mathrm{P}}$ is the matrix $[\boldsymbol{\rho}_1, \boldsymbol{\rho}_2, \boldsymbol{\rho}_3]$. After solving for $\epsilon$, we have our redistributed fluctuation waves,
\begin{align}
    A^+\Delta \tilde{Q}_{i,j} = \sum_{p: \omega_p > 0}\epsilon_p \omega_p \boldsymbol{\rho}_p \\
    A^-\Delta \tilde{Q}_{i,j} = \sum_{p: \omega_p< 0}\epsilon_p \omega_p \boldsymbol{\rho}_p.
\end{align}

Then we follow the algorithm in \cite{calhoun2008logically} and rotate them back:
\begin{align}
A^{\pm} \Delta Q_{i,j} = R_{i,j}^T A^{\pm}\Delta \tilde{Q}_{i,j}.
\end{align}
These waves are then weighted by the length of the barrier cut edge, as shown in \cref{fig:barst}.
This is then used to update the cut cells in a conservative but possibly unstable way,
\begin{align}
    Q^{L,n+1}_{i,j} & = Q^{L,n}_{i,j} - \frac{\Delta t}{\alpha^L_{i,j}}(l_{i,j}A^+ \Delta Q_{i,j} \nonumber \\
  & \hspace{-0.4cm} + l^L_{i-1/2,j}A^+\Delta Q^L_{i-1/2,j} + l^L_{i+1/2,j}A^-\Delta Q^L_{i+1/2,j} \nonumber \\&\hspace{-0.4cm} + l^L_{i,j+1/2} B^-\Delta Q^L_{i,j+1/2} + l^L_{i,j-1/2}B^+\Delta Q^L_{i,j-1/2}),
\end{align}
\begin{align}
Q^{U,n+1}_{i,j} & = Q^{U,n}_{i,j} - \frac{\Delta t}{\alpha^U_{i,j}}(l_{i,j}A^- \Delta Q_{i,j}  \nonumber \\
   & \hspace{-0.4cm}+ l^U_{i-1/2,j}A^+\Delta Q^U_{i-1/2,j} + l^U_{i+1/2,j}A^-\Delta Q^U_{i+1/2,j} \nonumber \\& \hspace{-0.4cm}+ l^U_{i,j+1/2} B^-\Delta Q^U_{i,j+1/2} + l^U_{i,j-1/2}B^+\Delta Q^U_{i,j-1/2}),
\end{align}
where $\alpha^{U/L}_{i,j}$ is the area of cut cell, $l_{i,j}$ the length of the barrier edge, and $l^{U/L}_{i\pm 1/2,j} $ and $l^{U/L}_{i,j\pm 1/2}$ represent the lengths of the vertical and horizontal edges of the cut cell, respectively. Note that there are five terms as the maximum number of cut cell edges is five, but some of them will drop depending on the shape of the cut cell (i.e. some $l$'s will be $0$).
\subsubsection{Neighborhood and overlap count}
Once we have updated each cut cell as described above, we then find the neighborhoods for each cut cell. In our model problems, finding the neighborhood for the cut cells is an easier task, as we can use normal neighborhoods for all of them.

\begin{figure}[h!]
    \centering
    \tikzset{every picture/.style={line width=0.75pt}} 
\begin{subfigure}[b]{0.4\textwidth}
\centering
\begin{tikzpicture}[scale=0.73,x=0.75pt,y=0.75pt,yscale=-1,xscale=1]
\draw (40,147.4) node [anchor=north west][inner sep=0.75pt]    {$j+1$};
\draw (61,197.4) node [anchor=north west][inner sep=0.75pt]    {$j$};

\draw (151,237.4) node [anchor=north west][inner sep=0.75pt]    {$i$};
\draw  [draw opacity=0] (81.6,81.6) -- (232.6,81.6) -- (232.6,232.6) -- (81.6,232.6) -- cycle ; \draw   (81.6,81.6) -- (81.6,232.6)(131.6,81.6) -- (131.6,232.6)(181.6,81.6) -- (181.6,232.6)(231.6,81.6) -- (231.6,232.6) ; \draw   (81.6,81.6) -- (232.6,81.6)(81.6,131.6) -- (232.6,131.6)(81.6,181.6) -- (232.6,181.6)(81.6,231.6) -- (232.6,231.6) ; \draw    ;
\draw  (231.6,109.6) -- (93.6,231.6); 
\draw [blue,fill=blue!10,opacity=0.8]  (181.6,153.6) -- (181.6,231.6) -- (131.6,231.6) -- (131.6,197.6) -- cycle ;
\draw [red,fill=red!10,opacity=0.8,thin] (131.6,197.6)--(181.6,153.6) -- (181.6,131.6) -- (131.6,131.6) -- (131.6,197.6) -- cycle ;
\end{tikzpicture}
\end{subfigure}
\begin{subfigure}[b]{0.4\textwidth}
\centering
\begin{tikzpicture}[x=0.75pt,y=0.75pt,yscale=-1,xscale=1]

\draw  [draw opacity=0] (188,55) -- (381,55) -- (381,151.5) -- (188,151.5) -- cycle ; \draw   (188,55) -- (188,151.5)(220,55) -- (220,151.5)(252,55) -- (252,151.5)(284,55) -- (284,151.5)(316,55) -- (316,151.5)(348,55) -- (348,151.5)(380,55) -- (380,151.5) ; \draw   (188,55) -- (381,55)(188,87) -- (381,87)(188,119) -- (381,119)(188,151) -- (381,151) ; \draw    ;
\draw    (188,79) -- (284,125.5) ;
\draw    (284,125.5) -- (380,77.5) ;
\draw  [red,fill=red!10,opacity=0.8,thin] (284,87) -- (284,125.5) -- (284,125.5) -- (252,110.25) -- (252,87) -- cycle ;
\draw  [red,fill=red!10,opacity=0.8,thin] (316,87) -- (316.5,109.75) -- (316.5,109.75) -- (284,125.5) -- (284,87) -- cycle ;
\draw  [blue,fill=blue!10,opacity=0.8] (252,110.25) -- (284,125.5) -- (284,151) -- (252,151) -- (252,151) -- cycle ;
\draw  [blue,fill=blue!10,opacity=0.8] (316.5,109.75) -- (316.5,109.75) -- (316,151) -- (284,151) -- (284,125.5) -- cycle ;

\draw (261,160) node [anchor=north west][inner sep=0.75pt]   [align=left] {$i-1$};
\draw (297.5,161.5) node [anchor=north west][inner sep=0.75pt]   [align=left] {$i$};
\draw (156.5,97.5) node [anchor=north west][inner sep=0.75pt]   [align=left] {$j+1$};
\draw (166.5,127) node [anchor=north west][inner sep=0.75pt]   [align=left] {$j$};

\end{tikzpicture}
\end{subfigure}

    \caption{Normal neighborhoods sufficient for model problems. On the top, we show upper neighborhood (in red) for cell ($i,j$) and lower neighborhood (in blue) for cell $(i,j+1)$ and on the bottom, we show upper neighborhoods (red) for cells ($i-1,j$) and ($i,j$) and lower neighborhoods (blue) for cells ($i-1,j+1$) and ($i,j+1$).}
    \label{fig:2Dnhood_mod}
\end{figure}
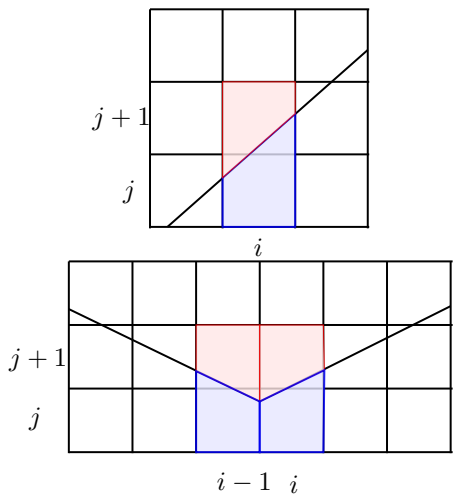

As shown in \cref{fig:2Dnhood_mod}, for each small cell, we choose the cell directly above or below as this will suffice to produce a neighborhood whose area is greater than $0.5\Delta x \Delta y$. For every cut cell produced by either the linear or V-barrier, there will always be a cell above (for upper cut cell) or below (for lower cut cell) whose $\alpha \ge 0.5$. This ensures using normal neighborhood for both the upper and lower half of the cut to be valid. Note that the obtuse angle of incidence for the V-barrier enables us to use normal neighborhoods extensively. Finding the neighborhood for the V-barrier essentially becomes a double copy of the linear problem. Note also that we must only take one sided neighborhoods to avoid mixing states across barrier.

The overlap count also becomes simplified in both model problems when we use the normal neighborhoods. For small cut cells on either side, we have an overlap count of $1$, as no other cut cell uses them for neighborhood formation. For the normal neighboring cells above or below small cut cells, we have an overlap count of $2$, as they have $\alpha \ge 0.5$ and are their own neighborhoods, in addition to being a neighbor to the small cut cells. If a cut cell has $\alpha \ge 0.5$, then it will automatically be its own neighborhood and have overlap count of $1$. Thus the overlap counts alternate between $1$ and $2$, depending on whether the cell is a small cut cell, non-small cut cell, or a neighbor of a small cut cell.

These overlap counts and the areas of the neighborhood cells then form neighborhood averages:
\begin{align}
    N^{U/L}_{i,j} = (\frac{1}{\alpha^{U/L}_{i,j} + \frac{\alpha^{U/L}_{i,j\pm 1}}{2\Delta x\Delta y}}) \big (\alpha^{U/L}_{i,j} Q^{U/L}_{i,j} + \frac{\alpha^{U/L}_{i,j\pm 1}}{2\Delta x\Delta y} Q^{U/L}_{i,j\pm 1} \big ),
\end{align}
for small cut cells $(i,j)$, where $Q^{U/L}_{i,j}$ represent the unstable but conservative update (the superscript $n+1$ is omitted for clarity and reserved for final update). For all other cells, we have $N^{U/L}_{i,j} = Q^{U/L}_{i,j}.$

All in all, we have for our update formula for all small cut cells ($i,j$):
\begin{align}
    Q^{n+1,U/L}_{i,j} = \begin{cases} N^{U/L}_{i,j} & \mbox{if } \text{overlap count} = 1\\
      \frac{1}{2} ( N^{U/L}_{i,j} + N^{U/L}_{i,j\pm 1} ) & \mbox{if } \text{overlap count} = 2 \label{2DSRDupd}
\end{cases}.
\end{align}

\section{Numerical examples}  \label{sec4}
All the numerical examples presented below have steady water height of $h=1.2$ and a dam jump of $\Delta h = 0.8$ or $1.5$, giving the overall height of the dam break to be $h_d=2.0$ or $2.7$. The barrier heights $\beta$ are chosen to be either $\beta=1.5$ or $\beta=5.0$. These are chosen to test both (1) complete reflection of incoming wave by the barrier ($\Delta h = 1.5 , \, \beta=5.0$) and (2) overtopping of wave over the barrier ($\Delta h = 0.8, \, \beta=1.5$). The boundary conditions for the problems are wall boundaries everywhere except for extrapolation condition on the overtopped side in the overtopping examples.

For comparison, we run two simulations. First is with the same initial condition in GeoClaw, a highly accurate geophysical fluid simulator \cite{BERGER20111195}, except with the barriers being single cell wide bathymetric jumps. We note that GeoClaw results are compared to ensure physical reliability of our results and that only similarity in quality of the numerical solution is sought after. Our problem is the numerical limit that the barrier thickness approaches zero, but it is difficult to attain reliable numerical solution using GeoClaw by thinning out the bathymetric jump. Therefore, we compare our simulations against the second, true numerical solutions performed in the mapped grid examples, where we use grid transformations and the wave redistribution method on a computational grid edge mapped to the barrier on the physical grid. We compare the 2D contour plots of water height at specific times from both the SRD cut cell method and GeoClaw/mapped examples and the gauge data, which are height measurements at a specified location in the domain (asterisked), in order to perform convergence analysis.

\subsection{The $20^\circ$ angled barrier}
\begin{figure}[h!]
\begin{subfigure}[b]{0.5\textwidth}
    \centering
    \includegraphics[scale=0.33]{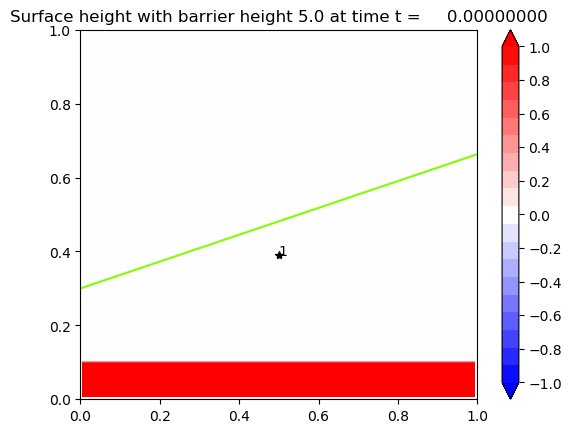}
    \caption{Reflection only case: $\beta = 5.0, h_d = 2.7$.}
    \label{fig:srf0}
\end{subfigure}
\begin{subfigure}[b]{0.5\textwidth}
     \centering
    \includegraphics[scale=0.33]{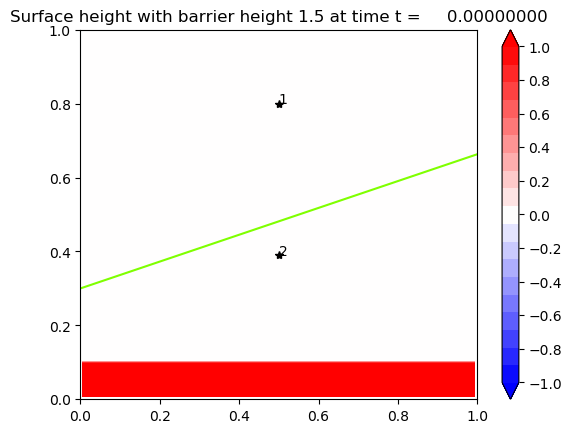}
    \caption{Overtopping case: $\beta = 1.5, h_d=2.0$.}
    \label{fig:sot0}
\end{subfigure}
    \caption{Initial condition for two problems. Both on 150$\times$150 grid.}
    \label{fig:Lprob_init}
\end{figure}

Here we present the first numerical example with a positively sloped 20$^\circ$ barrier. We first show a sufficiently high barrier example that prohibits overtopping of an incoming wave and then a lower barrier that permits overtopping (\cref{fig:Lprob_init}). In the reflection example, we place our gauge at $(0.5,0.39)$ just off the barrier, to capture the reflection close to the wall. In the overtopping example, we place two gauges, one at $(0.5,0.39)$ as before and another at $(0.5,0.8)$ to both capture the reflection and also the overtopped wave further away from the barrier.

\subsubsection{Case of Reflection Only}
We test a dam break problem to simulate the incidence of oblique waves upon the barrier.

\begin{figure}[h!]
\begin{subfigure}[b]{0.5\textwidth}

    \centering
    \includegraphics[scale=0.33]{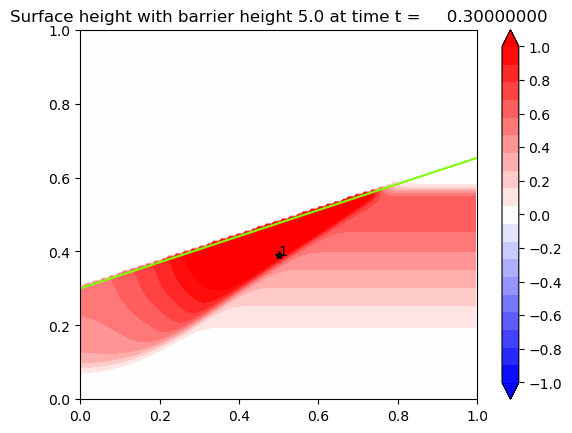}
    \caption{SRD}
    \label{fig:srf7}
    \end{subfigure}
    \begin{subfigure}[b]{0.5\textwidth}

    \centering
    \includegraphics[scale=0.33]{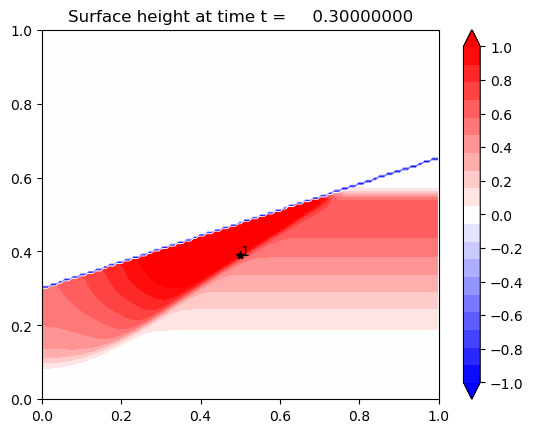}
    \caption{GeoClaw}
    \label{fig:srf7_gc}
    \end{subfigure}
      \caption{Comparison with GeoClaw at time $t=0.3$}
        \label{fig:srf7_fig}
\end{figure}

\begin{figure}[h!]
\begin{subfigure}[b]{0.5\textwidth}

    \centering
    \includegraphics[scale=0.33]{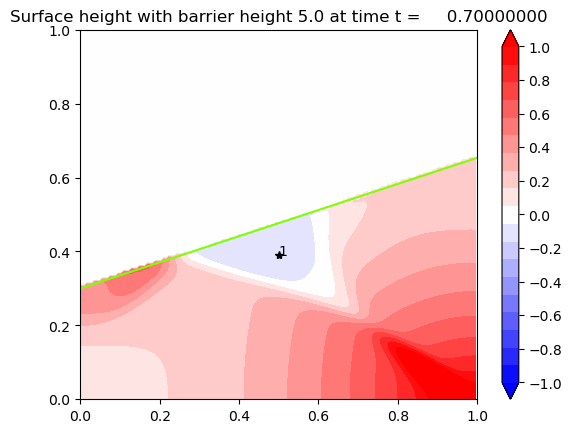}
    \caption{SRD}
    \label{fig:srf14}
    \end{subfigure}
    \begin{subfigure}[b]{0.5\textwidth}

    \centering
    \includegraphics[scale=0.33]{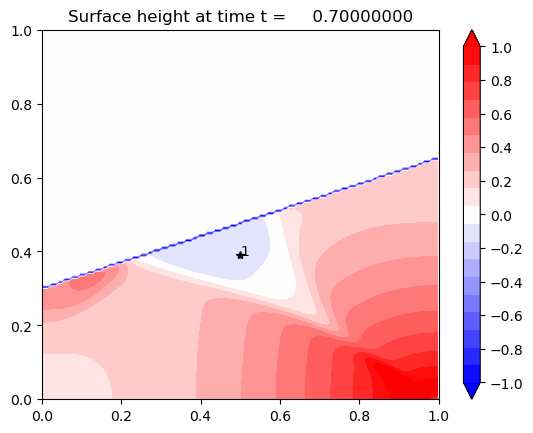}
    \caption{GeoClaw}
    \label{fig:srf14_gc}
    \end{subfigure}
      \caption{Comparison with GeoClaw at time $t=0.7.$}
        \label{fig:srf14_fig}
\end{figure}

In order to highlight the wave reflections we amplify the dam break size as well to $1.5.$ We observe very good comparison at both time snapshots (\cref{fig:srf7_fig} and \cref{fig:srf14_fig}).

\begin{figure}[h!]
    \centering
    \includegraphics[scale=0.35]{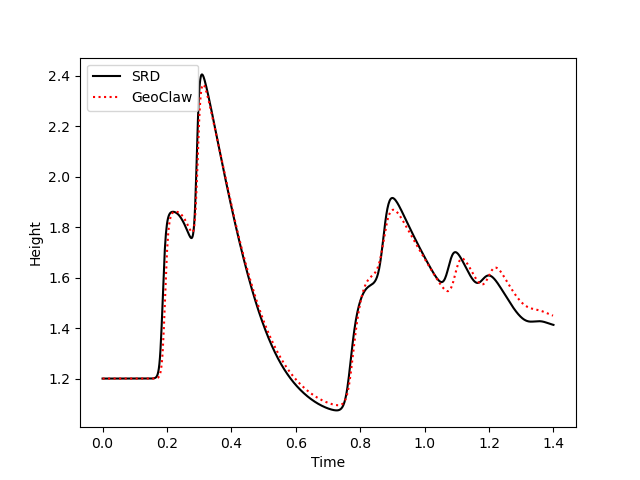}
    \caption{Time profile at Gauge 1 $(0.5,0.39)$ for complete blockage example with slanted barrier.}
    \label{fig:tprf}
\end{figure}

At time $t=0.3$, the incoming wave is gliding up along the barrier while being reflected back in direction normal to the barrier. One can observe how the wave front on the right has not yet reached the barrier. By time $t=0.7$, the reflected wave has bounced off the wall boundary condition on the bottom and reflected back onto the barrier, repeating the gliding motion as seen on the lower left region of the barrier. We can see that the gauge heights for all time are also in good agreement between the two methods (\cref{fig:tprf}).

\subsubsection{Case of Overtopping}
Next we test overtopping using the same slanted barrier. This time we take $\Delta h = 0.8$ and allow the the top edge of the domain to have extrapolation boundary.

\begin{figure}[h!]
\begin{subfigure}[b]{0.5\textwidth}
\centering
    \includegraphics[scale=0.33]{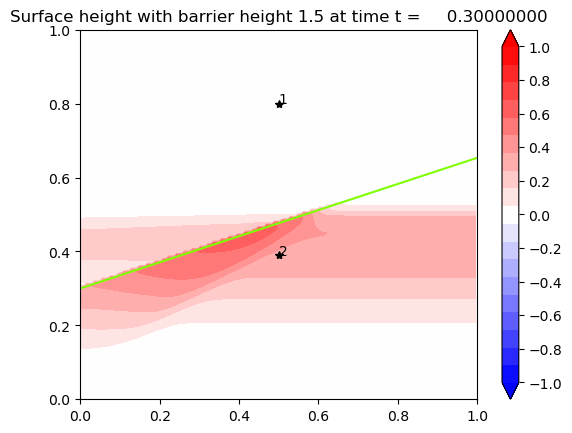}
    \caption{$t=0.3$}
    \label{fig:sot2}
\end{subfigure}
\begin{subfigure}[b]{0.5\textwidth}
\centering
    \includegraphics[scale=0.33]{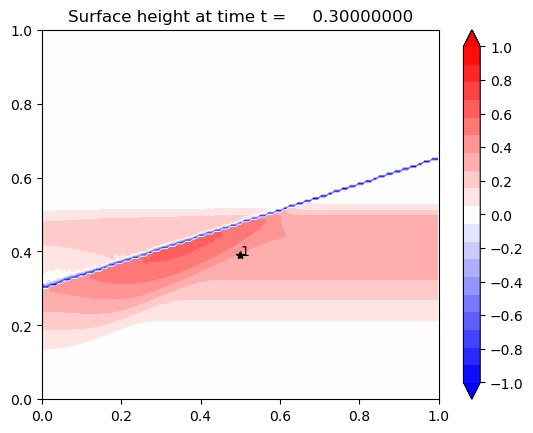}
    \caption{$t=0.3$}
    \label{fig:sot2gc}
\end{subfigure}
\caption{Comparison with GeoClaw at $t=0.3.$}
        \label{fig:sot2_fig}

\end{figure}

\begin{figure}[h!]
\begin{subfigure}[b]{0.5\textwidth}

    \centering
    \includegraphics[scale=0.33]{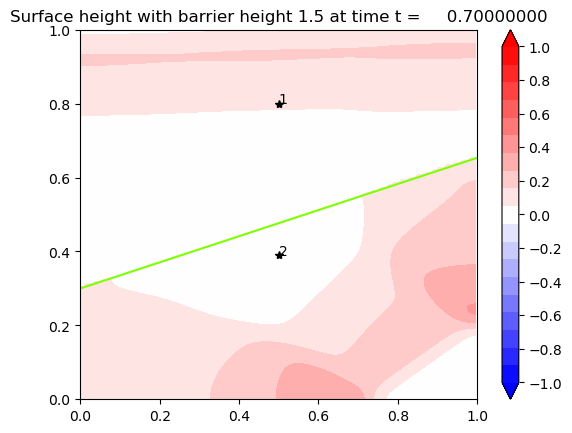}
    \caption{$t=0.7$}
    \label{fig:sot7}
    \end{subfigure}
    \begin{subfigure}[b]{0.5\textwidth}
    \centering
    \includegraphics[scale=0.33]{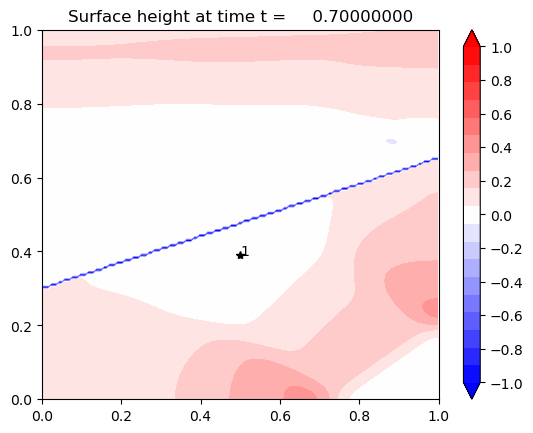}
    \caption{$t=0.7$}
    \label{fig:sot7_gc}
    \end{subfigure}
    \caption{Comparison with GeoClaw at $t=0.7.$}

        \label{fig:sot7_fig}

\end{figure}

We observe from both SRD results and GeoClaw results that the wave is abated from the barrier and proceeds in the same direction after overtopping (\cref{fig:sot2_fig}). The reflection shown at time $t=0.7$ on the lower side of the barrier is similar to the reflective behavior observed already in the complete blockage example (\cref{fig:sot7_fig}). An interesting observation to be made is at $t=1.4,$ the gliding waves `pinch up' to achieve overtopping momentum on the right, whereas it does not overtop on the left part of the barrier.

\begin{figure}[h!]
\begin{subfigure}[b]{0.5\textwidth}

    \centering
    \includegraphics[scale=0.33]{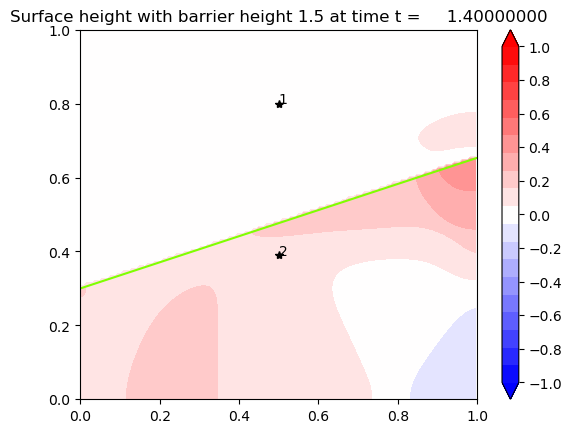}
    \caption{$t=1.4$}
    \label{fig:sot14}
    \end{subfigure}
    \begin{subfigure}[b]{0.5\textwidth}
    \centering
    \includegraphics[scale=0.33]{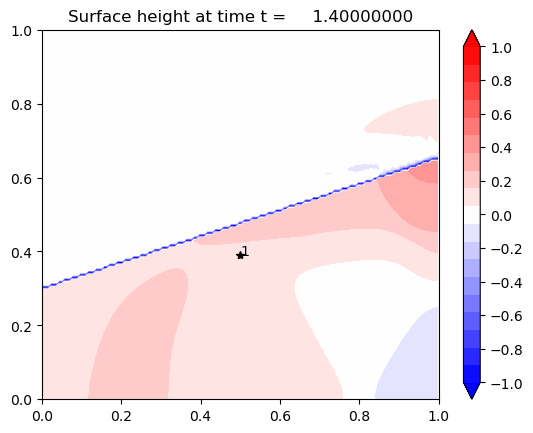}
    \caption{$t=1.4$}
    \label{fig:sot14_gc}
    \end{subfigure}
    \caption{Comparison with GeoClaw at $t=1.4.$}

        \label{fig:sot14_fig}

\end{figure}

\begin{figure}[h!]
\begin{subfigure}[b]{0.5\textwidth}
\centering
    \includegraphics[scale=0.34]{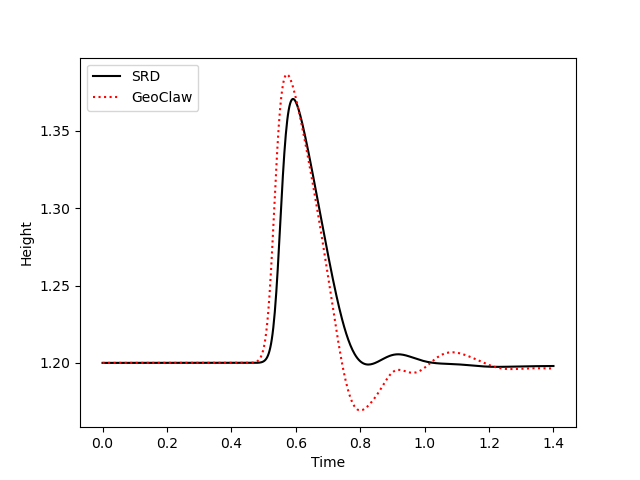}
    \caption{Time profile of Gauge 1 (0.5,0.8).}
    \label{fig:sottp}
\end{subfigure}
\begin{subfigure}[b]{0.5\textwidth}
\centering
    \includegraphics[scale=0.34]{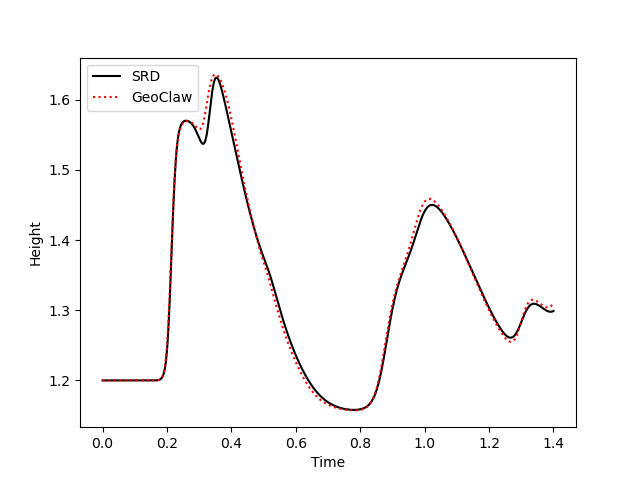}
    \caption{Time profile of Gauge 2 (0.5,0.39).}
    \label{fig:sottp2}
\end{subfigure}
\caption{Gauge profiles compared with GeoClaw results}
    \label{fig:sotTP}

\end{figure}
The similarity in behavior (\cref{fig:sot2_fig}-\cref{fig:sot14_fig}) and also magnitude of waves as seen in \cref{fig:sotTP} assure us of the physicality of the SRD results. We note that there is characteristic difference in the overtopped wave profile \cref{fig:sottp}. The GeoClaw example has a lower dip at the end of the first wave than the SRD. We attribute this to the presence of the bathymetric cell jump in the GeoClaw example causing a rarefaction at the ``behind'' interface of the jump (i.e. on the side of the overtopping) due to the large disparity in the heights \cite{george2008augmented}. Also the peaks are slightly off-centered due to the nonzero thickness of barrier causing earlier overtopping for GeoClaw. The reflected wave shown in \cref{fig:sottp2}, however, is virtually identical.

\subsubsection{Comparison to mapped grid}
\begin{figure}[h!]
    \centering
    \includegraphics[scale=0.35]{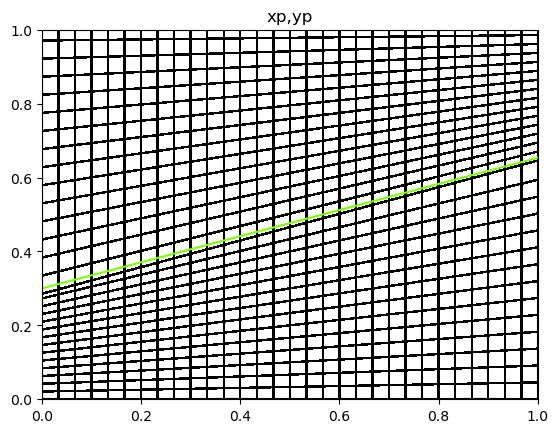}
    \caption{Mapped grid for the $20^\circ$ barrier. Coarsened to $30 \times 30$ to highlight mapping.}
    \label{fig:mapgridL}
\end{figure}

We validate the SRD results against a mapped grid example. We transform the computational uniform grid $(x,y)$ into a skewed grid $f_L(x,y)$ (\cref{fig:mapgridL}):
\begin{align}
    f_L(x,y) & = (x, \mu_L (y)), \\
    \mu_L(y) & =\begin{cases} \frac{L(x)}{y^*} y  &\mbox{if } y \in [0,y^*] \\
\frac{1-L(x)}{1-y^*}(y-1) +1 & \mbox{if } y \in [y^*,1]
 \end{cases},
\end{align}
where $L(x)$ is the barrier line equation and $y^*$ is the computational $y$-edge mapped to the barrier edge (lime green in \cref{fig:mapgridL}), with $y^* = \frac{y_1+y_2}{2}$ (\cref{fig:2d_setup1}). At $y=y^*$ wave redistribution is applied.

\begin{figure}[h!]
\begin{subfigure}[b]{0.5\textwidth}
\centering
    \includegraphics[scale=0.35]{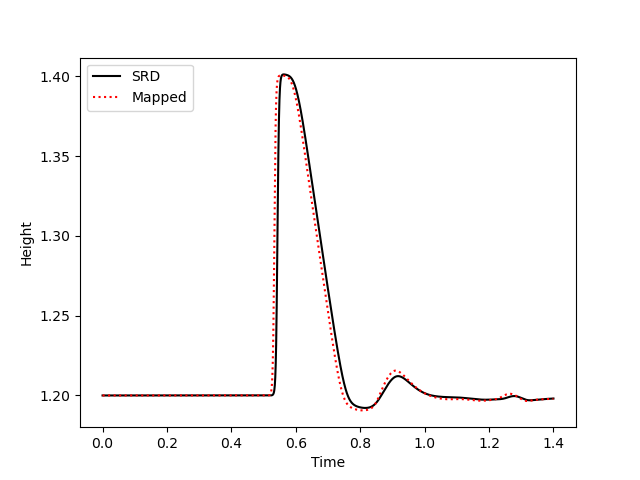}
    \caption{Time profile of Gauge 1 $(0.5,0.8)$.}
    \label{fig:otmp}
\end{subfigure}
\begin{subfigure}[b]{0.5\textwidth}
\centering
    \includegraphics[scale=0.35]{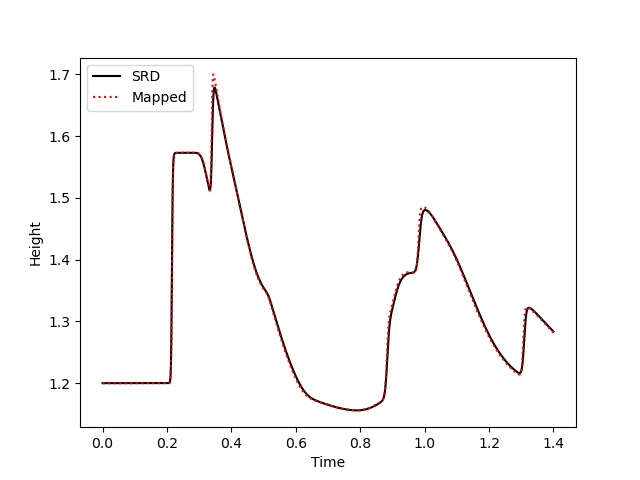}
    \caption{Time profile of Gauge 2 $(0.5,0.39)$.}
    \label{fig:otmp2}
\end{subfigure}
\caption{Gauge comparisons between SRD and mapped grid. Both on $900 \times 900$ grid.}
    \label{fig:OTMP}
\end{figure}

\begin{figure}[h!]
\begin{subfigure}[b]{0.5\textwidth}
\centering
    \includegraphics[scale=0.33]{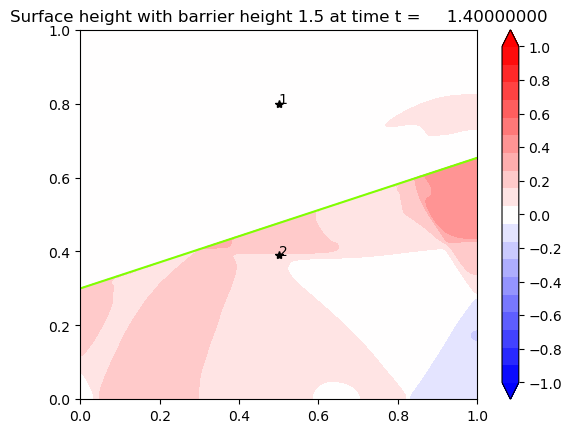}
    \caption{SRD results at $t=1.4$.}
    \label{fig:otsrd14}
\end{subfigure}
\begin{subfigure}[b]{0.5\textwidth}
\centering
    \includegraphics[scale=0.33]{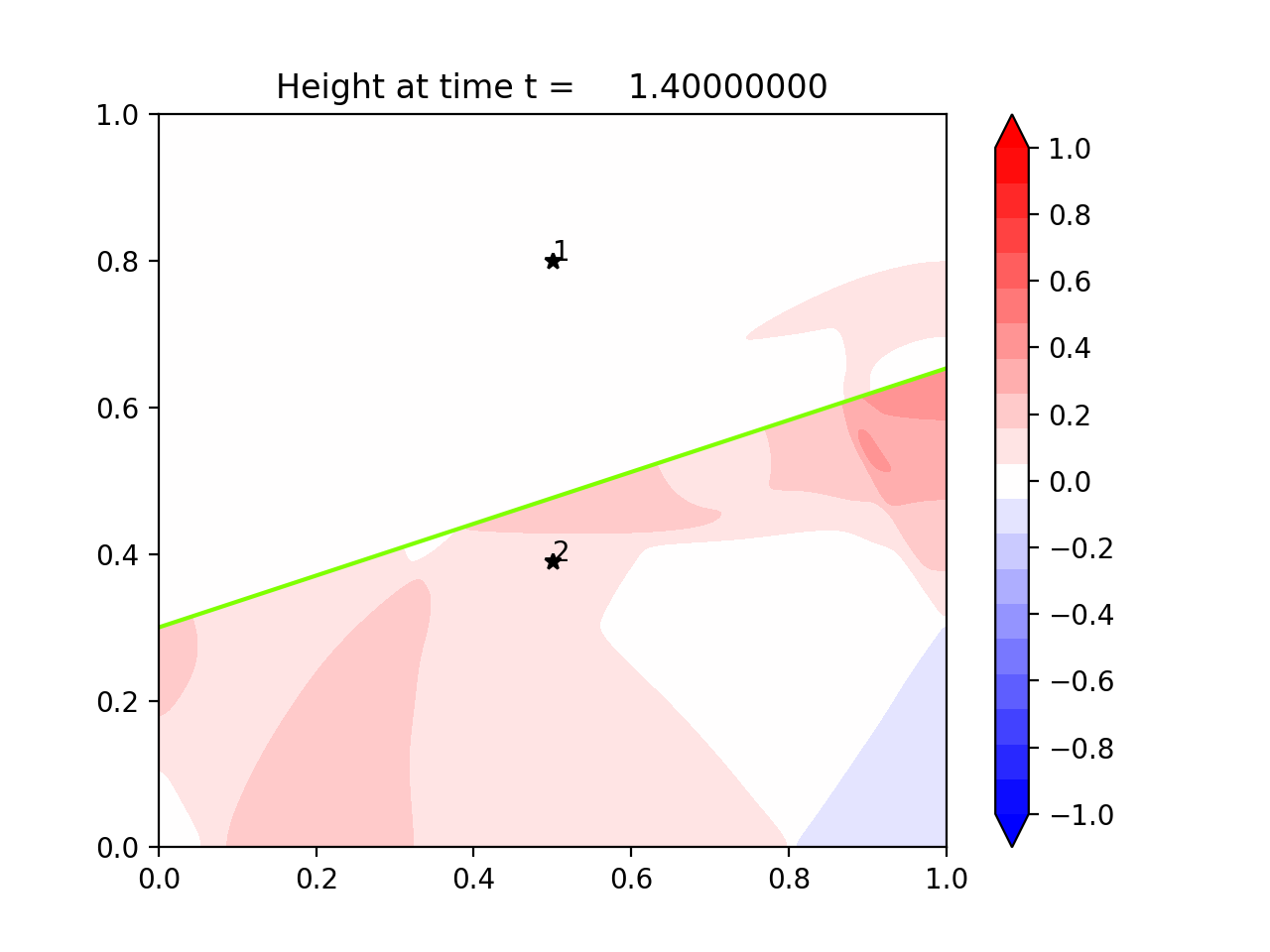}
    \caption{Mapped grid results at $t=1.4$.}
    \label{fig:otmp14}
\end{subfigure}
\caption{Comparison between SRD ($900 \times 900$) and mapped grid ($900 \times 900$) at $t=1.4$.}
\label{fig:srd_map_14}
\end{figure}

In \cref{fig:OTMP,fig:srd_map_14}, we observe the gauge results and contours plot at $t=1.4$ between the SRD and the mapped results both on $900\times 900$ grid and see the disappearance of the earlier differences.

\subsubsection{Convergence}
For convergence studies we observe the results at gauge point against the mapped grid ($900\times 900$) result, as would be of interest in storm simulations.
\begin{table}[h!]
\centering
\begin{tabular}{|c|c|c|c|}
\hline
~ & ~ & \multicolumn{2}{c|}{$L_1$ Error} \\
$\Delta x$ & $N_x, N_y$ & Gauge 1 & Gauge 2 \\ \hline
4.e-2   & 25    & 8.58e-3 & 2.70e-2\\ \hline
2.e-2   & 50    & 3.04e-3 & 1.14e-2\\ \hline
1.e-2   & 100   & 8.46e-4 & 3.98e-3\\ \hline
0.666e-2 & 150 & 3.89e-4 & 2.10e-3\\ \hline
0.333e-2 & 300 & 1.08e-4 & 5.23e-4\\ \hline
0.222e-2 & 450 & 5.37e-5 & 2.08e-4\\ \hline
0.1666-2 & 600 & 4.63e-5 & 1.01e-4\\ \hline
0.1333-2 & 750 & 3.44e-5 & 5.59e-5\\ \hline
0.1111-2 & 900 & 3.21e-5 & 3.14e-5\\ \hline

\end{tabular}
\caption{$L_1$ errors computed at gauge point 1 (0.5,0.8) and 2 (0.5, 0.39).}
\label{tab:Lbar_g1}
\end{table}
Also we only study the convergence of the overtopped examples, for similar reasons as above mentioned, namely, that in realistic scenarios the barriers will be overtopped by incoming waves.

\begin{figure}[h!]
    \begin{subfigure}[b]{0.5\textwidth}
    \centering
    \includegraphics[scale=0.4]{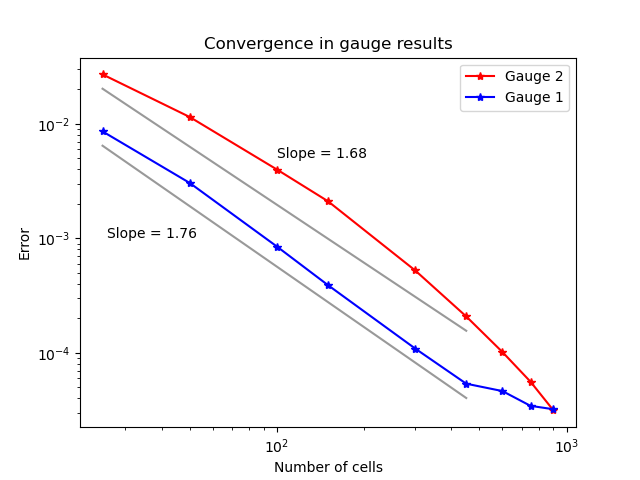}
    \caption{Against mapped grid: order around 1.7.}
    \label{fig:Lconv}
    \end{subfigure}
    \begin{subfigure}[b]{0.5\textwidth}
    \centering
    \includegraphics[scale=0.4]{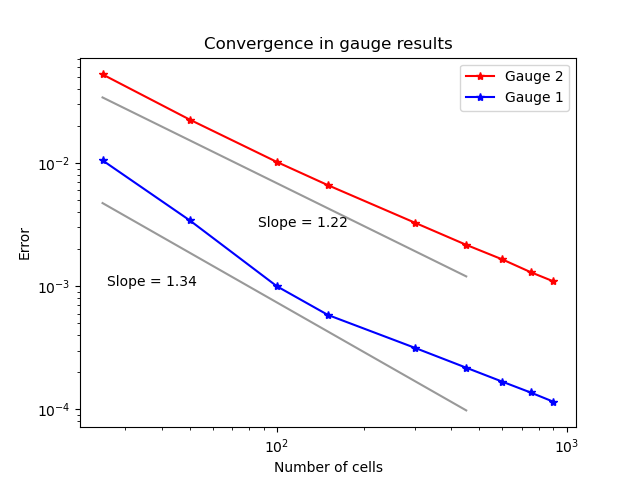}
    \caption{Against GeoClaw: order around 1.2-1.3.}
    \label{fig:Lconv_GC}
    \end{subfigure}
    \caption{Convergence plots with mapped grid (left) and GeoClaw results (right) for $20^\circ$ problem.}
    \label{fig:convL}
\end{figure}

The convergence is shown in both \cref{tab:Lbar_g1} and \cref{fig:Lconv}. We observe a convergence order of approximately 1.7 for both gauge points. We note that gauge 2 is closer to the barrier and yet the order of convergence is very similar to that for gauge 1, which is further away from the barrier.  This order of convergence is surprising, given that we are not doing any gradient reconstruction at the cut cells but only using piecewise constant values \cite{smith2007comparison}. This is most likely due to use of transverse solvers away from the barrier \cite{leveque2002finite}.

We also do convergence studies with the GeoClaw example by comparing our SRD results against GeoClaw example on $1200 \times 1200$ grid and the solution converges with greater than one order of convergence as seen in \cref{fig:Lconv_GC}.


\subsection{$117^\circ$ angled V-barrier}

\subsubsection{Case of Reflection Only}
\begin{figure}[h!]
\begin{subfigure}[b]{0.5\textwidth}
    \centering
    \includegraphics[scale=0.33]{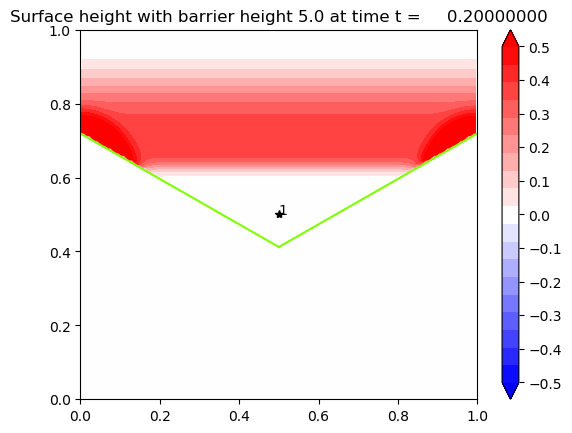}
    \caption{SRD: $t=0.2$}
    \label{fig:rf2}
    \end{subfigure}
    \begin{subfigure}[b]{0.5\textwidth}
      \centering
    \includegraphics[scale=0.33]{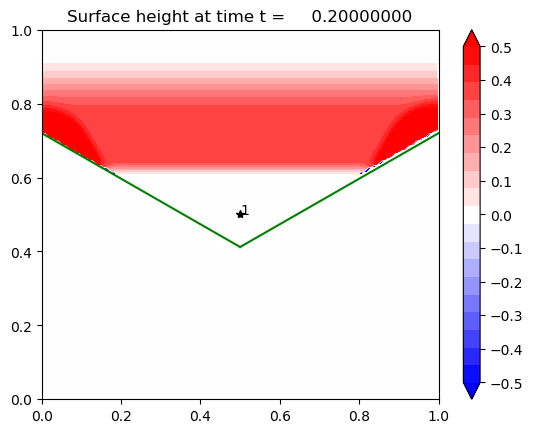}
    \caption{GeoClaw: $t=0.2$}
    \label{fig:rf2_gc}
    \end{subfigure}
    \caption{Because the barrier does not allow overtopping, the two problems are numerical the same.}
        \label{fig:rf2_fig}

\end{figure}

\begin{figure}[h!]
\begin{subfigure}[b]{0.5\textwidth}

    \centering
    \includegraphics[scale=0.33]{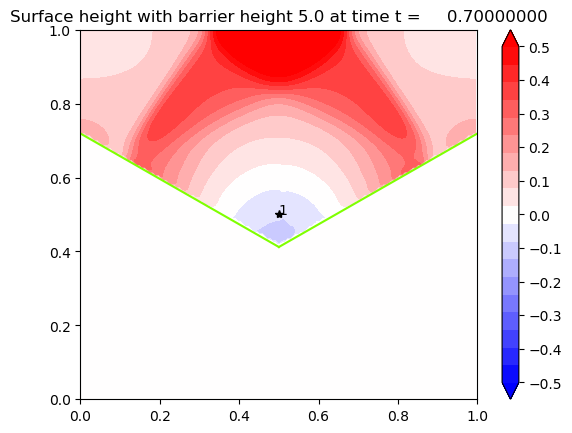}
    \caption{SRD: $t=0.7$}
    \label{fig:rf7}
    \end{subfigure}
    \begin{subfigure}[b]{0.5\textwidth}

    \centering
    \includegraphics[scale=0.33]{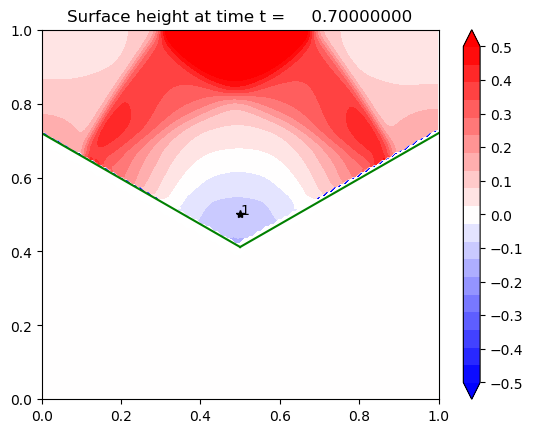}
    \caption{GeoClaw: $t=0.7$}
    \label{fig:rf7_gc}
    \end{subfigure}
      \caption{The gliding reflected waves have crossed each other in opposite directions.}
        \label{fig:rf7_fig}
\end{figure}

\begin{figure}[h!]
    \centering
    \includegraphics[scale=0.35]{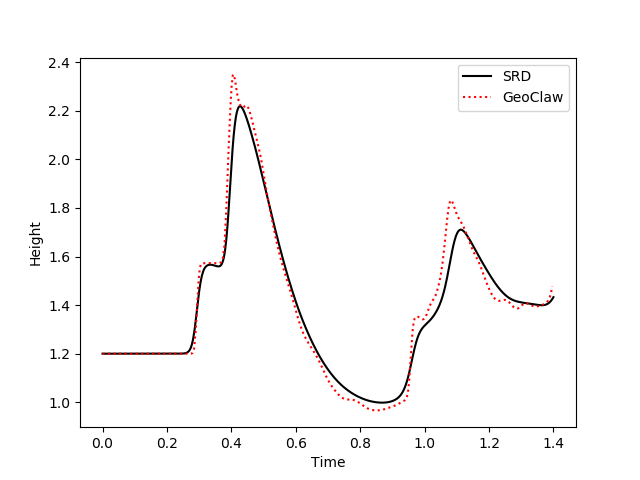}
    \caption{Time profile at Gauge 1 (0.5,0.5).}
    \label{fig:tprfV}
\end{figure}

\cref{fig:rf2_fig} to \cref{fig:rf7_fig} show an example with $150 \times 150$ grid in a domain of $[0,1] \times [0,1]$ for both SRD and GeoClaw and show a nice containment of the water behind the barrier. At \cref{fig:rf2_fig}, we can see the linear wave gliding across each side of the V-barrier towards the center. At \cref{fig:rf7_fig}, the gliding waves have crossed each other and are spreading radially outward away from the center of the barrier, leaving a dip (in blue). In \cref{fig:tprfV}, we show the wave profile at gauge point 1. We observe two waves: first initial wave that is reflected at the center and second wave reflected from the side boundary conditions. The second wave is lower in amplitude as it has lost some of its momentum.


\subsubsection{Case of Overtopping}
\begin{figure}[h!]
\begin{subfigure}[b]{0.5\textwidth}
\centering
    \includegraphics[scale=0.35]{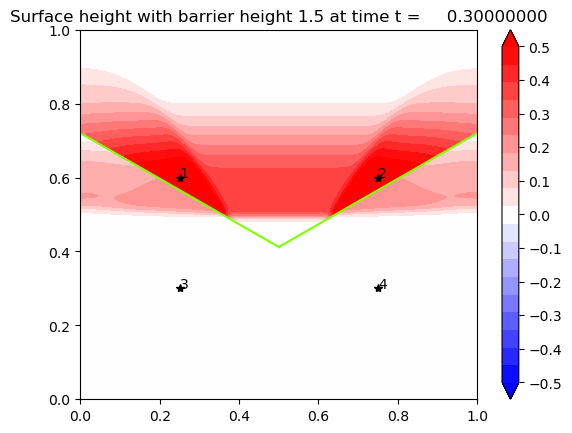}
    \caption{SRD: $t=0.3$}
    \label{fig:ot2}
\end{subfigure}
\begin{subfigure}[b]{0.5\textwidth}
\centering
    \includegraphics[scale=0.35]{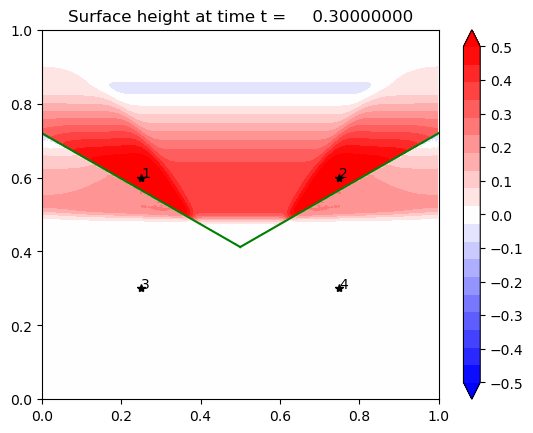}
    \caption{GeoClaw: $t=0.3$}
    \label{fig:ot2gc}
\end{subfigure}
\caption{Both $300 \times 300$ grid. Note the structure of the just overtopped wave.}
        \label{fig:ot2_fig}

\end{figure}

\begin{figure}[h!]
\begin{subfigure}[b]{0.5\textwidth}

    \centering
    \includegraphics[scale=0.35]{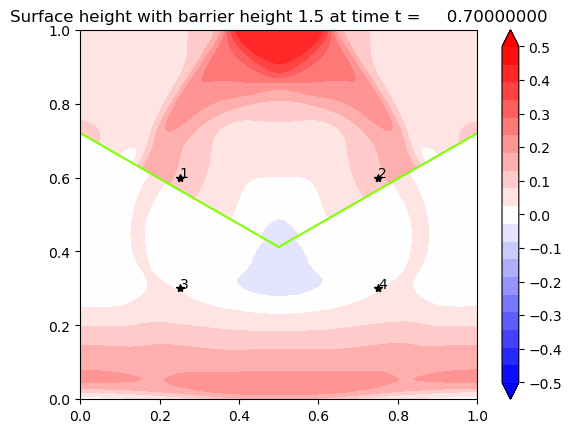}
    \caption{SRD: $t=0.7$}
    \label{fig:ot7}
    \end{subfigure}
    \begin{subfigure}[b]{0.5\textwidth}
    \centering
    \includegraphics[scale=0.35]{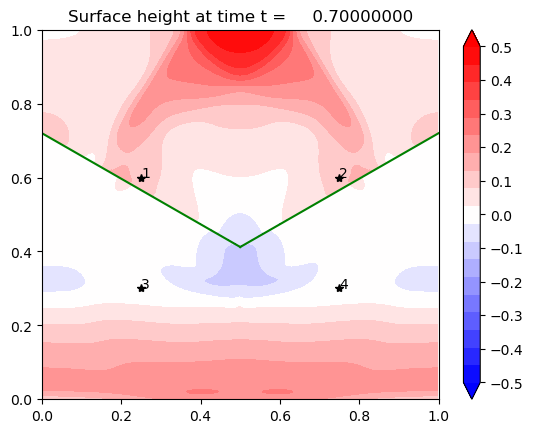}
    \caption{GeoClaw: $t=0.7$}
    \label{fig:ot7_gc}
    \end{subfigure}
    \caption{Note the radially outward moving wave from the center of the V-barrier.}

        \label{fig:ot7_fig}

\end{figure}

\begin{figure}[h!]
\begin{subfigure}[b]{0.5\textwidth}

    \centering
    \includegraphics[scale=0.35]{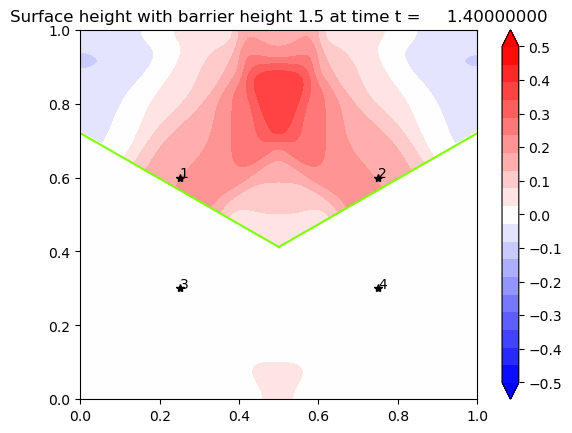}
    \caption{SRD: $t=1.4$}
    \label{fig:ot14}
    \end{subfigure}
    \begin{subfigure}[b]{0.5\textwidth}
    \centering
    \includegraphics[scale=0.35]{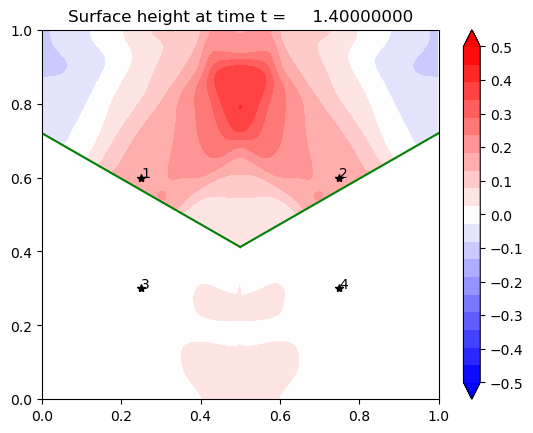}
    \caption{GeoClaw: $t=1.4$}
    \label{fig:ot14_gc}
    \end{subfigure}
    \caption{Note the ``island'' of peak at the bottom center.}

        \label{fig:ot14_fig}

\end{figure}

We increase the resolution to 300 for clearer results and first do a comparison against GeoClaw results to show physicality of our results. From \cref{fig:ot2_fig} to \cref{fig:ot14_fig}, we can see the qualitative similarities of the two 2D plots. At time $t=0.3$, we see the overtopping wave's ``wing''-like structure just below the V-barrier where the amplitude is highest. In \cref{fig:ot7_fig}, we see the overtopping wave moving radially outward from the center, more captured in the SRD results than GeoClaw. Finally in \cref{fig:ot14_fig} we see the small island of wave amplitude at the bottom center of the plot in both results.

\begin{figure}[h!]
\begin{subfigure}[b]{0.5\textwidth}
\centering
    \includegraphics[scale=0.35]{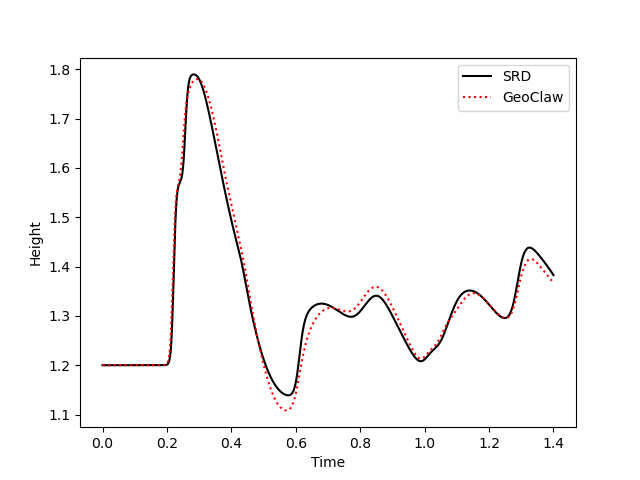}
    \caption{Time profile of Gauge 1,2 $(0.25(0.75),0.6)$.}
    \label{fig:ottp}
\end{subfigure}

\begin{subfigure}[b]{0.5\textwidth}
\centering
    \includegraphics[scale=0.35]{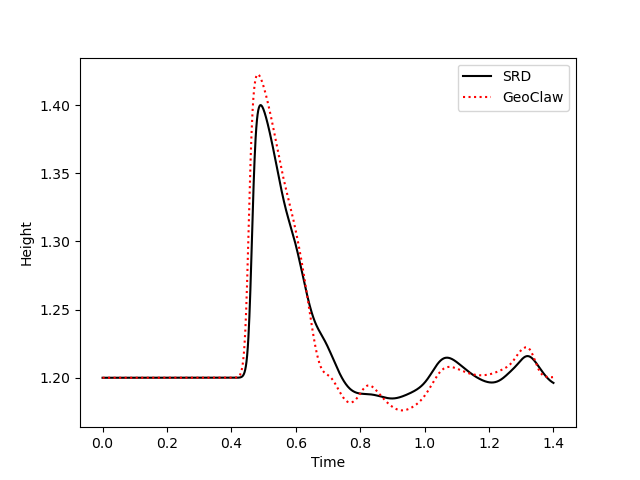}
    \caption{Time profile of Gauge 3,4 $(0.25(0.75),0.3)$.}
    \label{fig:ottp3}
\end{subfigure}
\caption{Gauge profiles compared with GeoClaw results: $300 \times 300$ grid.}
    \label{fig:otTP}
\end{figure}

We place our gauges at either side of the V-barrier $(0.25,0.3)$, $(0.75,0.3)$, $(0.25,0.6)$, $(0.75,0.6)$ in order to test for symmetry in the results. Indeed we do find symmetry as can be seen in the idential plots of the gauge results in \cref{fig:otTP}. We observe the major peaks lining up between the two results. Again we see slightly earlier peak and deeper dip in the overtopped waves, due to the presence of physical barrier.

\subsubsection{Computational Superiority to Refinement using GeoClaw}
\label{sec5}
To highlight the superiority of using the SRD method on the zero width barrier, we do a V-barrier simulation on GeoClaw using same resolution as SRD but adaptive (double) refinement at the barrier. In reality the refinement level required at barriers will be greater than 2 as barriers are much skinnier than surrounding bathymetric surfaces. The barrier in these GeoClaw runs is two cells wide, to get down to single cell width in the adaptive mesh refinement.

We show results from using two resolutions $\Delta x = 1/300, \, 1/450$. These already show the computational benefit we derive from our proposed method. We compare the $\Delta t$'s and number of steps taken from both SRD and GeoClaw simulations. For $\Delta x =1/300$ we observe that we get tenfold increase in the minimum $\Delta t$ (from 8.6e-06 to 9.6e-05) and 320 \% increase in the average $\Delta t$ and about fivefold decrease in the number of steps taken (from 0.00028 to 0.0012 and 9958 steps to 1850 steps). For $\Delta x = 1/450$ we observe about 70-fold increase in the minimum $\Delta t$ (from 2.9e-06 to 2.2e-04) and 340 \% increase in average $\Delta t$ and 5.5 times reduction in number of steps taken (from 1.7e-4 to 7.7e-4 and 15861 to 2843 steps).
\subsection{Comparison to mapped grid}
\begin{figure}[h!]
    \centering
    \includegraphics[scale=0.35]{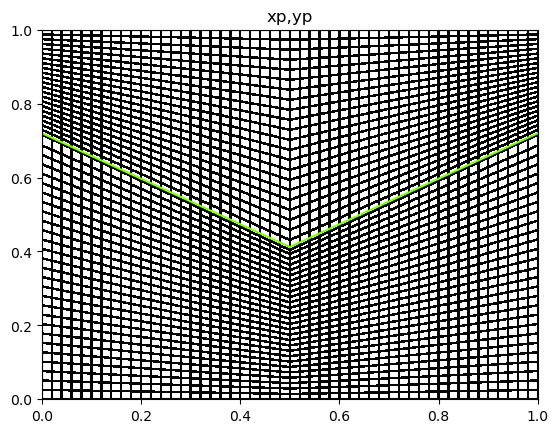}
    \caption{Mapped grid for the V barrier. Coarsened to $50 \times 50$ to highlight mapping.}
    \label{fig:mapgrid}
\end{figure}

The mapped grid is shown in \cref{fig:mapgrid}.
Here we transform the computational uniform grid into a physical V shaped grid akin to what is done in \cite{berger2012simplified}, with the following mapping $f$:
{\small
\begin{align}
    f_V(x,y) & = (x, \mu_L (y)), \\
    \mu_V(y) & =\begin{cases} \frac{L_1(x)}{y^*} y, & (x,y)\in [0,0.5]\times [0,y^*] \\
\frac{1-L_1(x)}{1-y^*}(y-1) +1,  & [0,0.5]\times [y^*,1] \\
\frac{L_2(x)}{y^*} y,&[0.5,1]\times [0,y^*] \\
\frac{1-L_2(x)}{1-y^*}(y-1) +1,& [0.5,1]\times [y^*,1]
 \end{cases},
\end{align}
}%
where $y^*$ is the computational barrier edge taken to be $y^* = 0.5(y_1+y_2)$ (\cref{fig:2d_setup2}), and $L_1(x)$ is the barrier line from coordinate 1 to 2 and $L_2(x)$ that from coordinate 2 to 3. We apply wave redistribution (\cref{WR}) at $y^*$.

\subsubsection{Convergence}
\begin{figure}[h!]
\begin{subfigure}[b]{0.5\textwidth}
\centering
    \includegraphics[scale=0.35]{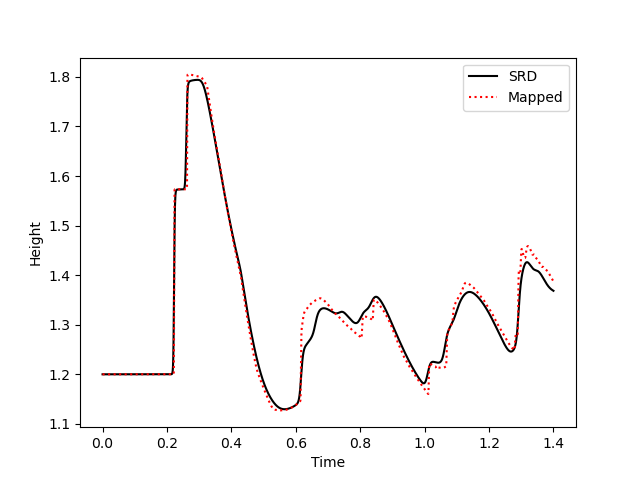}
    \caption{Time profile of Gauge 1,2 $(0.25(0.75),0.6)$.}
    \label{fig:ottpV}
\end{subfigure}

\begin{subfigure}[b]{0.5\textwidth}
\centering
    \includegraphics[scale=0.35]{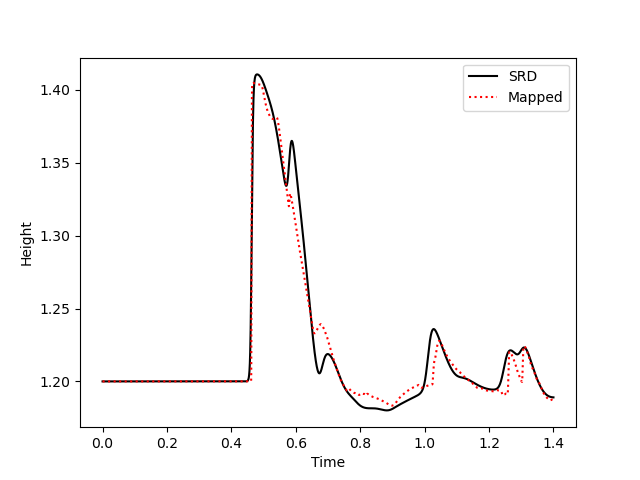}
    \caption{Time profile of Gauge 3,4 $(0.25(0.75),0.3)$.}
    \label{fig:ottpV3}
\end{subfigure}
\caption{Gauge profiles compared with mapped grid results: $1050 \times 1050$ for SRD and $1250 \times 1250$ for mapped grid.}
    \label{fig:otTPV}
\end{figure}
\begin{figure}[h!]
    \begin{subfigure}[b]{0.5\textwidth}
    \centering
    \includegraphics[scale=0.4]{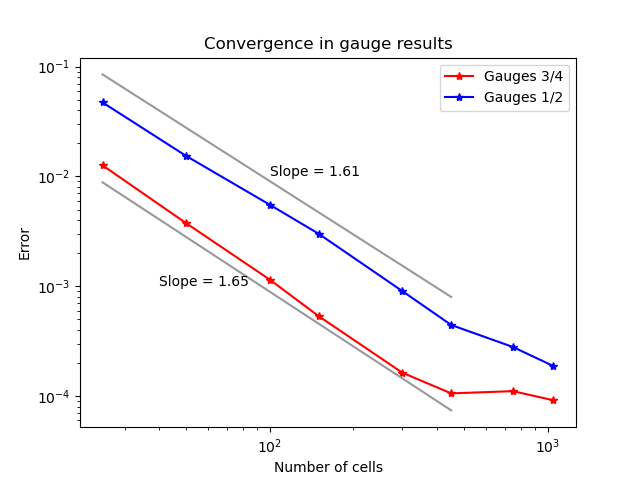}
    \caption{Against mapped grid: order around 1.6.}
    \label{fig:Vconv}
    \end{subfigure}
    \begin{subfigure}[b]{0.5\textwidth}
     \centering
    \includegraphics[scale=0.4]{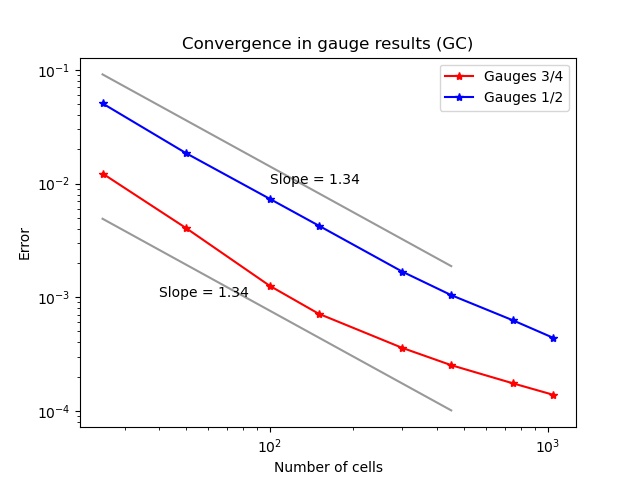}
    \caption{Against GeoClaw: order around 1.3.}
    \label{fig:VconvGC}
    \end{subfigure}
    \caption{Convergence plots with mapped grid (left) and GeoClaw results (right)}
    \label{fig:convV}
\end{figure}

In \cref{fig:otTPV} we plot the gauge results of $1250 \times 1250$ mapped grid V-barrier example and $1050 \times 1050$ SRD example. We do see that the SRD results contain more fine movements of the wave and that the mapped grid example produces more smooth wave patterns. Overall, however, we see convergence as shown in \cref{tab:Vbar_g1,fig:Vconv} (for gauge points 1 through 4). The order of convergence are somewhere around 1.6 for both gauge points (\cref{fig:Vconv}). Again we observe greater-than-one order of convergence despite using piecewise constant approximations.

\begin{figure}[h!]
\begin{subfigure}[b]{0.5\textwidth}
\centering
    \includegraphics[scale=0.35]{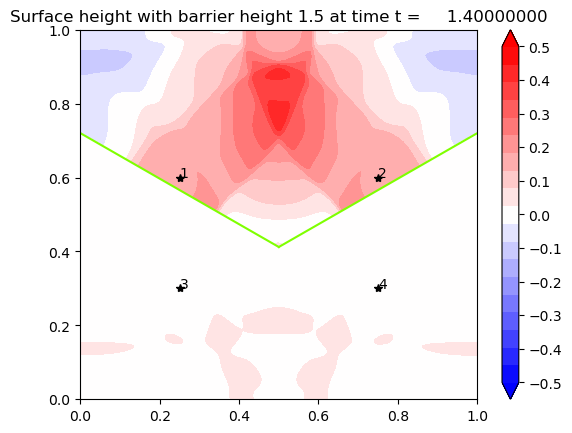}
    \caption{SRD results at $t=1.4$.}
    \label{fig:otVsrd14}
\end{subfigure}
\begin{subfigure}[b]{0.5\textwidth}
\centering
    \includegraphics[scale=0.35]{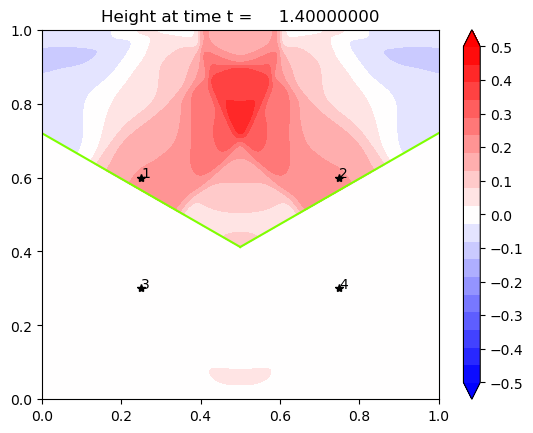}
    \caption{Mapped grid results at $t=1.4$.}
    \label{fig:otVmp14}
\end{subfigure}
\caption{Comparison between SRD ($1050 \times 1050$) and mapped grid ($1250 \times 1250$) results.}
    \label{fig:otV14}
\end{figure}

Finally, in \cref{fig:otV14} we show two 2D plots comparing the SRD results on $1050 \times 1050$ grid and mapped grid results on $1250 \times 1250$ grid at $t=1.4.$ We see that SRD is less diffusive even on a slightly lower resolution that the mapped grid example. This can be seen in the finer details in both the reflected side (where the reflected waves cross in the center) and also the overtopped side, in the appearance of more residual overtopped waves.

We also compare our V-barrier SRD results with the finite width thin barrier simulation on GeoClaw, refined to $1200 \times 1200$ grid, i.e. $\Delta x  \approx 0.0009$ thick barrier and observe a convergence order of 1.3 as shown in \cref{fig:VconvGC}.

\begin{table}[h!]
\centering
\begin{tabular}{|c|c|c|c|}
\hline
~ & ~ & \multicolumn{2}{c|}{$L_1$ Error} \\
$\Delta x$ & $N_x, N_y$ & Gauge 1,2 & Gauge 3,4\\ \hline
4.e-2 & 25 & 4.72e-2 & 1.26e-2 \\ \hline
2.e-2 & 50 & 1.53e-2 & 3.75e-3 \\ \hline
1.e-2 & 100 & 5.52e-3 & 1.14e-3\\ \hline
.666e-2 & 150 & 3.00e-3 & 5.34e-4 \\ \hline
.333e-2 & 300 & 9.04e-4 & 1.63e-4 \\ \hline
.222e-2 & 450 & 4.44e-4 & 1.05e-4 \\ \hline
.1333-2 & 750 & 2.80e-4 & 1.10e-4 \\ \hline
.1111-2 & 1050 & 1.87e-4 & 9.16e-5 \\ \hline

\end{tabular}
\caption{$L_1$ errors computed as difference between mapped and SRD at gauge points 1,2 located at (0.25,0.6) and (0.75,0.6) and points 3,4 located at (0.25,0.3) and (0.75,0.3).}
\label{tab:Vbar_g1}
\end{table}



\section{Conclusion}

We have developed an SRD method applicable for storm barrier modeling. Two model barriers have been given, one being a linear segment and the other being a double linear segment forming a V-shape. We observe both qualitative and quantitative similarities in the results of three methods: GeoClaw (using single cell wide barriers), mapped grid, and the SRD. Most importantly, in SRD we significantly increase the average $\Delta t$ and also the $\min \Delta t$ of the barrier simulations due to both the cut cell method and zero width approximation, saving computational time from refining the barrier. Furthermore, the method does not involve the complexities of the $h$-box method which calculates lot of geometrical parameters. However, we obtain an order greater than 1, despite using only piecewise constant values at every cell, without using any gradient reconstruction.

\backmatter

\bmhead{Code availability}

To access the underlying codes for the simulations, please refer to github.com/cr2940

\bibliography{sn-bibliography}


\begin{thebibliography}{25}
\ifx \bisbn   \undefined \def \bisbn  #1{ISBN #1}\fi
\ifx \binits  \undefined \def \binits#1{#1}\fi
\ifx \bauthor  \undefined \def \bauthor#1{#1}\fi
\ifx \batitle  \undefined \def \batitle#1{#1}\fi
\ifx \bjtitle  \undefined \def \bjtitle#1{#1}\fi
\ifx \bvolume  \undefined \def \bvolume#1{\textbf{#1}}\fi
\ifx \byear  \undefined \def \byear#1{#1}\fi
\ifx \bissue  \undefined \def \bissue#1{#1}\fi
\ifx \bfpage  \undefined \def \bfpage#1{#1}\fi
\ifx \blpage  \undefined \def \blpage #1{#1}\fi
\ifx \burl  \undefined \def \burl#1{\textsf{#1}}\fi
\ifx \doiurl  \undefined \def \doiurl#1{\url{https://doi.org/#1}}\fi
\ifx \betal  \undefined \def \betal{\textit{et al.}}\fi
\ifx \binstitute  \undefined \def \binstitute#1{#1}\fi
\ifx \binstitutionaled  \undefined \def \binstitutionaled#1{#1}\fi
\ifx \bctitle  \undefined \def \bctitle#1{#1}\fi
\ifx \beditor  \undefined \def \beditor#1{#1}\fi
\ifx \bpublisher  \undefined \def \bpublisher#1{#1}\fi
\ifx \bbtitle  \undefined \def \bbtitle#1{#1}\fi
\ifx \bedition  \undefined \def \bedition#1{#1}\fi
\ifx \bseriesno  \undefined \def \bseriesno#1{#1}\fi
\ifx \blocation  \undefined \def \blocation#1{#1}\fi
\ifx \bsertitle  \undefined \def \bsertitle#1{#1}\fi
\ifx \bsnm \undefined \def \bsnm#1{#1}\fi
\ifx \bsuffix \undefined \def \bsuffix#1{#1}\fi
\ifx \bparticle \undefined \def \bparticle#1{#1}\fi
\ifx \barticle \undefined \def \barticle#1{#1}\fi
\bibcommenthead
\ifx \bconfdate \undefined \def \bconfdate #1{#1}\fi
\ifx \botherref \undefined \def \botherref #1{#1}\fi
\ifx \url \undefined \def \url#1{\textsf{#1}}\fi
\ifx \bchapter \undefined \def \bchapter#1{#1}\fi
\ifx \bbook \undefined \def \bbook#1{#1}\fi
\ifx \bcomment \undefined \def \bcomment#1{#1}\fi
\ifx \oauthor \undefined \def \oauthor#1{#1}\fi
\ifx \citeauthoryear \undefined \def \citeauthoryear#1{#1}\fi
\ifx \endbibitem  \undefined \def \endbibitem {}\fi
\ifx \bconflocation  \undefined \def \bconflocation#1{#1}\fi
\ifx \arxivurl  \undefined \def \arxivurl#1{\textsf{#1}}\fi
\csname PreBibitemsHook\endcsname

\bibitem{usarmy}
\begin{botherref}
\oauthor{\bparticle{Corps~of} \bsnm{Engineers}, \binits{U.A.}}:
Hurricane sandy aftermath: Hurricane barriers managed by corps engineers in new
  england prevent 29.7 million in damages.
New England District
(2012)
\end{botherref}
\endbibitem

\bibitem{jmse8090725}
\begin{botherref}
\oauthor{\bsnm{Chen}, \binits{Z.}},
\oauthor{\bsnm{Orton}, \binits{P.}},
\oauthor{\bsnm{Wahl}, \binits{T.}}:
Storm surge barrier protection in an era of accelerating sea-level rise:
  Quantifying closure frequency, duration and trapped river flooding.
Journal of Marine Science and Engineering
\textbf{8}(9)
(2020).
\doiurl{10.3390/jmse8090725}
\end{botherref}
\endbibitem

\bibitem{BERGER20111195}
\begin{barticle}
\bauthor{\bsnm{Berger}, \binits{M.J.}},
\bauthor{\bsnm{George}, \binits{D.L.}},
\bauthor{\bsnm{LeVeque}, \binits{R.J.}},
\bauthor{\bsnm{Mandli}, \binits{K.T.}}:
\batitle{The GeoClaw software for depth-averaged flows with adaptive
  refinement}.
\bjtitle{Advances in Water Resources}
\bvolume{34}(\bissue{9}),
\bfpage{1195}--\blpage{1206}
(\byear{2011}).
\doiurl{10.1016/j.advwatres.2011.02.016}.
\bcomment{New Computational Methods and Software Tools}
\end{barticle}
\endbibitem

\bibitem{martin2010lake}
\begin{botherref}
\oauthor{\bsnm{Martin}, \binits{S.K.}},
\oauthor{\bsnm{Savant}, \binits{G.}},
\oauthor{\bsnm{McVan}, \binits{D.C.}}:
Lake borgne surge barrier study.
Technical report,
ENGINEER RESEARCH AND DEVELOPMENT CENTER VICKSBURG MS COASTAL AND HYDRAULICS
  LAB
(2010)
\end{botherref}
\endbibitem

\bibitem{zhang2013transition}
\begin{barticle}
\bauthor{\bsnm{Zhang}, \binits{K.}},
\bauthor{\bsnm{Li}, \binits{Y.}},
\bauthor{\bsnm{Liu}, \binits{H.}},
\bauthor{\bsnm{Rhome}, \binits{J.}},
\bauthor{\bsnm{Forbes}, \binits{C.}}:
\batitle{Transition of the coastal and estuarine storm tide model to an
  operational storm surge forecast model: A case study of the florida coast}.
\bjtitle{Weather and forecasting}
\bvolume{28}(\bissue{4}),
\bfpage{1019}--\blpage{1037}
(\byear{2013})
\end{barticle}
\endbibitem

\bibitem{berger2012simplified}
\begin{barticle}
\bauthor{\bsnm{Berger}, \binits{M.}},
\bauthor{\bsnm{Helzel}, \binits{C.}}:
\batitle{A simplified h-box method for embedded boundary grids}.
\bjtitle{SIAM Journal on Scientific Computing}
\bvolume{34}(\bissue{2}),
\bfpage{861}--\blpage{888}
(\byear{2012})
\end{barticle}
\endbibitem

\bibitem{2006-Chung-p607}
\begin{barticle}
\bauthor{\bsnm{Chung}, \binits{M.-H.}}:
\batitle{{C}artesian cut cell approach for simulating incompressible flows with
  rigid bodies of arbitrary shape}.
\bjtitle{computers \& fluids}
\bvolume{35},
\bfpage{607}--\blpage{623}
(\byear{2006}).
\doiurl{10.1016/j.compfluid.2005.04.005}
\end{barticle}
\endbibitem

\bibitem{colella2006cartesian}
\begin{barticle}
\bauthor{\bsnm{Colella}, \binits{P.}},
\bauthor{\bsnm{Graves}, \binits{D.T.}},
\bauthor{\bsnm{Keen}, \binits{B.J.}},
\bauthor{\bsnm{Modiano}, \binits{D.}}:
\batitle{A cartesian grid embedded boundary method for hyperbolic conservation
  laws}.
\bjtitle{Journal of Computational Physics}
\bvolume{211}(\bissue{1}),
\bfpage{347}--\blpage{366}
(\byear{2006})
\end{barticle}
\endbibitem

\bibitem{BERGER20171}
\begin{botherref}
Chapter 1 - cut cells: Meshes and solvers.
In: \oauthor{\bsnm{Abgrall}, \binits{R.}},
\oauthor{\bsnm{Shu}, \binits{C.-W.}} (eds.)
Handbook of Numerical Methods for Hyperbolic Problems.
Handbook of Numerical Analysis,
vol. 18,
pp. 1--22.
Elsevier
(2017)
\end{botherref}
\endbibitem

\bibitem{INGRAM2003561}
\begin{barticle}
\bauthor{\bsnm{Ingram}, \binits{D.M.}},
\bauthor{\bsnm{Causon}, \binits{D.M.}},
\bauthor{\bsnm{Mingham}, \binits{C.G.}}:
\batitle{Developments in cartesian cut cell methods}.
\bjtitle{Mathematics and Computers in Simulation}
\bvolume{61}(\bissue{3}),
\bfpage{561}--\blpage{572}
(\byear{2003}).
\doiurl{10.1016/S0378-4754(02)00107-6}.
\bcomment{MODELLING 2001 - Second IMACS Conference on Mathematical Modelling
  and Computational Methods in Mechanics, Physics, Biomechanics and
  Geodynamics}
\end{barticle}
\endbibitem

\bibitem{BERGER2020109820}
\begin{botherref}
\oauthor{\bsnm{Berger}, \binits{M.}},
\oauthor{\bsnm{Giuliani}, \binits{A.}}:
A state redistribution algorithm for finite volume schemes on cut cell meshes.
Journal of Computational Physics,
109820
(2020).
\doiurl{10.1016/j.jcp.2020.109820}
\end{botherref}
\endbibitem

\bibitem{calhoun2008logically}
\begin{barticle}
\bauthor{\bsnm{Calhoun}, \binits{D.A.}},
\bauthor{\bsnm{Helzel}, \binits{C.}},
\bauthor{\bsnm{LeVeque}, \binits{R.J.}}:
\batitle{Logically rectangular grids and finite volume methods for pdes in
  circular and spherical domains}.
\bjtitle{SIAM review}
\bvolume{50}(\bissue{4}),
\bfpage{723}--\blpage{752}
(\byear{2008})
\end{barticle}
\endbibitem

\bibitem{GDavid}
\begin{barticle}
\bauthor{\bsnm{Leveque}, \binits{R.}},
\bauthor{\bsnm{George}, \binits{D.}},
\bauthor{\bsnm{Berger}, \binits{M.}}:
\batitle{Tsunami modelling with adaptively refined finite volume methods}.
\bjtitle{Acta Numerica}
\bvolume{20},
\bfpage{211}--\blpage{289}
(\byear{2011}).
\doiurl{10.1017/S0962492911000043}
\end{barticle}
\endbibitem

\bibitem{leveque2002finite}
\begin{bbook}
\bauthor{\bsnm{LeVeque}, \binits{R.J.}}, \betal:
\bbtitle{Finite Volume Methods for Hyperbolic Problems},
(\byear{2002})
\end{bbook}
\endbibitem

\bibitem{CAUSON2000545}
\begin{barticle}
\bauthor{\bsnm{Causon}, \binits{D.M.}},
\bauthor{\bsnm{Ingram}, \binits{D.M.}},
\bauthor{\bsnm{Mingham}, \binits{C.G.}},
\bauthor{\bsnm{Yang}, \binits{G.}},
\bauthor{\bsnm{Pearson}, \binits{R.V.}}:
\batitle{Calculation of shallow water flows using a cartesian cut cell
  approach}.
\bjtitle{Advances in Water Resources}
\bvolume{23}(\bissue{5}),
\bfpage{545}--\blpage{562}
(\byear{2000}).
\doiurl{10.1016/S0309-1708(99)00036-6}
\end{barticle}
\endbibitem

\bibitem{TUCKER2000591}
\begin{barticle}
\bauthor{\bsnm{Tucker}, \binits{P.G.}},
\bauthor{\bsnm{Pan}, \binits{Z.}}:
\batitle{A cartesian cut cell method for incompressible viscous flow}.
\bjtitle{Applied Mathematical Modelling}
\bvolume{24}(\bissue{8}),
\bfpage{591}--\blpage{606}
(\byear{2000}).
\doiurl{10.1016/S0307-904X(00)00005-6}
\end{barticle}
\endbibitem

\bibitem{may2017explicit}
\begin{barticle}
\bauthor{\bsnm{May}, \binits{S.}},
\bauthor{\bsnm{Berger}, \binits{M.}}:
\batitle{An explicit implicit scheme for cut cells in embedded boundary
  meshes}.
\bjtitle{Journal of Scientific Computing}
\bvolume{71}(\bissue{3}),
\bfpage{919}--\blpage{943}
(\byear{2017})
\end{barticle}
\endbibitem

\bibitem{bale2003wave}
\begin{barticle}
\bauthor{\bsnm{Bale}, \binits{D.S.}},
\bauthor{\bsnm{Leveque}, \binits{R.J.}},
\bauthor{\bsnm{Mitran}, \binits{S.}},
\bauthor{\bsnm{Rossmanith}, \binits{J.A.}}:
\batitle{A wave propagation method for conservation laws and balance laws with
  spatially varying flux functions}.
\bjtitle{SIAM Journal on Scientific Computing}
\bvolume{24}(\bissue{3}),
\bfpage{955}--\blpage{978}
(\byear{2003})
\end{barticle}
\endbibitem

\bibitem{einfeldt1988godunov}
\begin{barticle}
\bauthor{\bsnm{Einfeldt}, \binits{B.}}:
\batitle{On godunov-type methods for gas dynamics}.
\bjtitle{SIAM Journal on Numerical Analysis}
\bvolume{25}(\bissue{2}),
\bfpage{294}--\blpage{318}
(\byear{1988})
\end{barticle}
\endbibitem

\bibitem{roe1981approximate}
\begin{barticle}
\bauthor{\bsnm{Roe}, \binits{P.L.}}:
\batitle{Approximate riemann solvers, parameter vectors, and difference
  schemes}.
\bjtitle{Journal of computational physics}
\bvolume{43}(\bissue{2}),
\bfpage{357}--\blpage{372}
(\byear{1981})
\end{barticle}
\endbibitem

\bibitem{causon2000calculation}
\begin{barticle}
\bauthor{\bsnm{Causon}, \binits{D.M.}},
\bauthor{\bsnm{Ingram}, \binits{D.M.}},
\bauthor{\bsnm{Mingham}, \binits{C.G.}},
\bauthor{\bsnm{Yang}, \binits{G.}},
\bauthor{\bsnm{Pearson}, \binits{R.V.}}:
\batitle{Calculation of shallow water flows using a cartesian cut cell
  approach}.
\bjtitle{Advances in water resources}
\bvolume{23}(\bissue{5}),
\bfpage{545}--\blpage{562}
(\byear{2000})
\end{barticle}
\endbibitem

\bibitem{BERGER2015180}
\begin{barticle}
\bauthor{\bsnm{Berger}, \binits{M.}}:
\batitle{A note on the stability of cut cells and cell merging}.
\bjtitle{Applied Numerical Mathematics}
\bvolume{96},
\bfpage{180}--\blpage{186}
(\byear{2015}).
\doiurl{10.1016/j.apnum.2015.05.003}
\end{barticle}
\endbibitem

\bibitem{li2021h}
\begin{barticle}
\bauthor{\bsnm{Li}, \binits{J.}},
\bauthor{\bsnm{Mandli}, \binits{K.T.}}:
\batitle{An h-box method for shallow water equations including barriers}.
\bjtitle{SIAM Journal on Scientific Computing}
\bvolume{43}(\bissue{2}),
\bfpage{431}--\blpage{454}
(\byear{2021})
\end{barticle}
\endbibitem

\bibitem{george2008augmented}
\begin{barticle}
\bauthor{\bsnm{George}, \binits{D.L.}}:
\batitle{Augmented riemann solvers for the shallow water equations over
  variable topography with steady states and inundation}.
\bjtitle{Journal of Computational Physics}
\bvolume{227}(\bissue{6}),
\bfpage{3089}--\blpage{3113}
(\byear{2008})
\end{barticle}
\endbibitem

\bibitem{smith2007comparison}
\begin{bchapter}
\bauthor{\bsnm{Smith}, \binits{T.}},
\bauthor{\bsnm{Barone}, \binits{M.}},
\bauthor{\bsnm{Bond}, \binits{R.}},
\bauthor{\bsnm{Lorber}, \binits{A.}},
\bauthor{\bsnm{Baur}, \binits{D.}}:
\bctitle{Comparison of reconstruction techniques for unstructured mesh vertex
  centered finite volume schemes}.
In: \bbtitle{18th AIAA Computational Fluid Dynamics Conference},
p. \bfpage{3958}
(\byear{2007})
\end{bchapter}
\endbibitem

\end{thebibliography}


\end{document}